\documentclass{article}

\usepackage{etex}
\usepackage{graphicx}
\usepackage{latexsym}
\usepackage{amsfonts}
\usepackage{amssymb}
\usepackage{amsmath}
\usepackage{verbatim}
\usepackage{url}
\usepackage[standard, thmmarks]{ntheorem}
\usepackage[margin=3cm]{geometry}
\usepackage{fancyhdr}

\newtheorem{thm}{Theorem}[section]

\newtheorem{lem}[thm]{Lemma}
\newtheorem{prop}[thm]{Proposition}

\newtheorem{qn}[thm]{Question}

\newcommand{\R}{\mathbb R}

\newcommand{\Z}{\mathbb Z}
\newcommand{\C}{\mathbb C}

\newcommand{\F}{\mathcal F}

\newcommand{\To}{\longrightarrow}

\newcommand{\E}{\mathcal{E}}
\newcommand{\z}{{\bf z}}
\newcommand{\B}{\mathcal{B}}
\newcommand{\D}{\mathcal{D}}
\newcommand{\Pent}{\mathcal{P}}
\newcommand{\V}{\mathcal{V}}
\newcommand{\T}{\mathcal{T}}

\newcommand{\hyp}{\mathbb H}
\newcommand{\Out}{\text{Out}\,}

\DeclareMathOperator{\Par}{Par} \DeclareMathOperator{\Hyp}{Hyp}
\DeclareMathOperator{\Ell}{Ell} 

 \DeclareMathOperator{\mcg}{MCG}
\DeclareMathOperator{\Aut}{Aut} \DeclareMathOperator{\Inn}{Inn}
 \DeclareMathOperator{\Isom}{Isom}
\DeclareMathOperator{\Axis}{Axis} \DeclareMathOperator{\Twist}{Tw}
\DeclareMathOperator{\Tr}{Tr} 
\DeclareMathOperator{\Fix}{Fix}
\DeclareMathOperator{\Ad}{Ad}

\begin{document}

\title{Hyperbolic cone-manifold structures with prescribed holonomy II: higher genus}

\author{Daniel V. Mathews}%

\date{}

\maketitle
% ----------------------------------------------------------------

\begin{abstract}

We consider the relationship between hyperbolic cone-manifold structures on surfaces, and algebraic representations of the fundamental group into a group of isometries. A hyperbolic cone-manifold structure on a surface, with all interior cone angles being integer multiples of $2\pi$, determines a holonomy representation of the fundamental group. We ask, conversely, when a representation of the fundamental group is the holonomy of a hyperbolic cone-manifold structure. In this paper we build upon previous work with punctured tori to prove results for higher genus surfaces.

Our techniques construct fundamental domains for hyperbolic cone-manifold structures, from the geometry of a representation. Central to these techniques are the Euler class of a representation, the group $\widetilde{PSL_2\R}$, the twist of hyperbolic isometries, and character varieties. We consider the action of the outer automorphism and related groups on the character variety, which is measure-preserving with respect to a natural measure derived from its symplectic structure, and ergodic in certain regions. Under various hypotheses, we almost surely or surely obtain a hyperbolic cone-manifold structure with prescribed holonomy.

\end{abstract}

\tableofcontents

\section{Introduction}

\subsection{The story continues...}

In this series of papers we consider geometric structures and holonomy representations. A geometric $(X, \Isom X)$ structure on an orientable manifold $M$ induces a holonomy representation $\rho: \pi_1 (M) \To \Isom X$; we ask, conversely, given a representation $\rho: \pi_1(M) \To \Isom X$, is $\rho$ the holonomy of a geometric structure? We consider 2-dimensional hyperbolic geometry and allow cone singularities; a representation $\rho$ then only makes sense if every interior cone point has a cone angle which is an integer multiple of $2\pi$.

This paper is a continuation of \cite{Me10MScPaper1}: see that paper for background and context for this problem. There we showed precisely which homomorphisms $\rho: \pi_1(S) \To PSL_2\R$ are holonomy representations, in the above sense, for $S$ a punctured torus. Precisely, we proved that $\rho$ is the holonomy representation of a hyperbolic cone-manifold structure on $S$ with geodesic boundary, except for at most one corner point, and no interior cone points, if and only if $\rho$ is not virtually abelian.

In this paper we extend this result, applying the ideas of \cite{Me10MScPaper0} and \cite{Me10MScPaper1} to surfaces other than the punctured torus, in particular higher genus surfaces.  Our results are not as complete as in the punctured torus case but may still be of interest. We will make use of the notion of twist of a hyperbolic isometry, studied in \cite{Me10MScPaper0} and which has properties relating it to areas in the hyperbolic plane and to the algebraically-defined angle function $\Theta$ of a matrix developed by Milnor \cite{Milnor}. This will give us nice relationships between representations of the fundamental group, their Euler class, and the arrangements of isometries in the hyperbolic plane encoded by them. We will make use of the main theorem of \cite{Me10MScPaper1} to construct hyperbolic cone-manifold structures on punctured tori: the idea is that, given a more complicated surface, we can cut it into punctured tori and pairs of pants, construct hyperbolic structures on the pieces, and glue them together. For this, we will need a similar result constructing hyperbolic structures on pairs of pants: this is section \ref{sec:hyperbolize_pants}.

This paper also proves a result that \emph{almost all} representations of a certain type are holonomy representations. The proof uses certain \emph{ergodicity} properties of the action of the mapping class group on the character variety \cite{Goldman97, Goldman03}. As such this gives a geometric application of results that previously may have been solely of dynamic and analytic interest; as far as we know it is the first such application.

As such, this paper carries out several tasks. It will show you how to hyperbolize your pants. It will give some background on the Euler class, character varieties, and actions on them. And it will show, when circumstances are favourable, with probability $1$ or with certainty, how to cut a surface into pieces, and piece by piece, hyperbolize it all.

\subsection{Results and discussion}

The answer to the present question --- which representations of a surface group into $PSL_2\R$ are holonomy representations? --- in the absence of cone points, i.e. for complete hyperbolic structures with totally geodesic (or cusped) boundary, is a theorem of Goldman \cite{Goldman_thesis}. For $S$ a closed surface with $\chi(S) < 0$, a representation $\rho: \pi_1(S) \To PSL_2\R$ determines an \emph{Euler class} $\E(\rho)$. (We discuss the Euler class in more detail below in sections \ref{sec:euler_class} and \ref{sec:algebraic_description_of_euler_class}.) By the Milnor--Wood inequality \cite{Milnor, Wood, Me10MScPaper0}, the Euler class $\E(\rho)$, evaluated on the fundamental class $[S]$, is no more than $\chi(S)$ in magnitude. The Euler class parametrises the connected components of the $PSL_2\R$-representation space \cite{Goldman88}. Goldman in \cite{Goldman_thesis} proved that $\rho$ is the holonomy of a hyperbolic structure on $S$ if and only if the Euler class is ``extremal'', i.e. $\pm \chi(S)$ times the fundamental cohomology class. Thus holonomy representations form precisely two of the $2|\chi(S)|+1$ components of the representation space. If $S$ has boundary, then the same machinery applies, and the same theorem holds, provided that each boundary curve is sent to a non-elliptic isometry. In this case we obtain a \emph{relative} Euler class. 

In this paper we will, \emph{inter alia}, reprove Goldman's theorem using our own methods, which are quite different from Goldman's. As it turns out, assuming an extremal Euler class forces the arrangement of isometries in the holonomy group to be highly favourable to our constructions; we can then obtain fundamental domains for punctured tori and pairs of pants, which glue together. In this sense our proof is perhaps more low-powered than Goldman's: we ``construct fundamental domains by hand''. However we must rely on a theorem of Gallo--Kapovich--Marden in \cite{Gallo_Kapovich_Marden} guaranteeing decompositions along hyperbolic curves, for which the proof is long and detailed.
\begin{thm}[Goldman \cite{Goldman_thesis}]
\label{goldman_theorem}
Let $S$ be a compact connected orientable surface with $\chi(S) < 0$, and let $\rho$ be a homomorphism $\pi_1(S) \To PSL_2\R$. If $S$ has boundary, assume $\rho$ takes each boundary curve to a non-elliptic element, so the relative Euler class $\E(\rho)$ is well-defined. The following are equivalent:
\begin{enumerate}
\item
$\rho$ is the holonomy of a complete hyperbolic structure on $S$ with totally geodesic or cusped boundary components (respectively as each boundary curve is taken by $\rho$ to a hyperbolic or parabolic);
\item
$\E(\rho)[S] = \pm \chi(S)$.
\end{enumerate}
\end{thm}
In the case of a closed surface $S$ of genus $g \geq 2$, Goldman's theorem simply becomes that a representation $\rho: \pi_1(S) \To PSL_2\R$ is the holonomy of a complete hyperbolic structure on $S$ if and only if $\E(\rho)[S] = \pm \chi(S)$.

As an aside, note that Goldman's theorem involves lifting $\rho$ to the universal cover of the isometry group. Reviewing the background in \cite{Me10MScPaper1}, lifts to universal covers also feature in corresponding results for other geometries. The Euler class can be considered as an \emph{obstruction} to lifting $\rho$ into the universal cover $\widetilde{PSL_2\R}$; see section \ref{sec:algebraic_description_of_euler_class} below. So from Goldman's result, a holonomy representation does not lift to $\widetilde{PSL_2\R}$; in fact, such representations are as ``un-liftable" as possible. This is in contrast to other geometries.

When cone points are introduced, the problem becomes more subtle. However there is a simple relation between that the Euler class of the holonomy representation and the number and type of cone points: the defect from extremeness of the Euler class must be made up by adding extra cone points (proposition \ref{euler_class_of_cone_manifold}). But among the components of the space of $PSL_2\R$ representations with non-extremal Euler class, it is not yet clear which among them are holonomy representations. So far as we know, it is still an open question whether the set of holonomy representations is dense among representations of Euler class $\pm 1, \pm 2, \ldots, \pm (\chi(S) + 1)$. Using ergodicity methods, as we do in this paper, one might hope that the set of such representations is conull, i.e. almost every representation is a holonomy representation in these cases.

We can however obtain some results as follows. When we can apply our punctured torus theorem \cite{Me10MScPaper1} twice, ``back to back'', we obtain a theorem proving the existence of hyperbolic cone-manifold structures on the genus 2 closed surface. 
\begin{thm}
\label{restricted_version}
    Let $S$ be a genus 2 closed surface. Let $\rho: \pi_1(S) \To PSL_2\R$ be a representation with
    $\E(\rho)[S] = \pm 1$. Suppose that there is a separating curve
    $C$ on $S$ such that $\rho(C)$ is not hyperbolic. Then $\rho$ is
    the holonomy of a hyperbolic cone-manifold structure on $S$ with
    one cone point of angle $4\pi$.
\end{thm}

And when circumstances permit us to apply the ergodicity results of \cite{Goldman03}, we have the following.
\begin{thm}
\label{technical_prop}
    Let $S$ be a closed orientable surface of genus $g \geq 2$.
    Almost every representation $\rho:
    \pi_1(S) \To PSL_2\R$ with $\E(\rho)[S] = \pm (\chi(S)+1)$,
    which sends some non-separating simple closed curve $C$
    to an elliptic, is the holonomy of a
    hyperbolic cone-manifold structure on $S$ with a single cone point with
    cone angle $4\pi$.
\end{thm}
As we proceed, we will introduce a measure on the character variety, so that this statement makes sense: the precise statement is theorem \ref{thm:technical_precise}.

Our proof relies on the existence of a non-separating simple closed curve with elliptic holonomy; this can then be cut off, allowing us to localise the deficiency in the Euler class. But we have not been able to show such a curve exists in general. Perhaps the ``almost" can be removed as well; despite a comment by Tan \cite{Tan94}, as far as we know the question remains open.
\begin{qn}
    For a general surfaces, is almost every representation with Euler class $\pm (\chi(S) + 1)$ a holonomy representation? Every representation?
\end{qn}

Our approach relies heavily on $\E(\rho)$ being close to extremal: once we localise the deficiency in $\E(\rho)$, we cut it off, and the rest of the representation has extremal Euler class. For other values of $\E(\rho)$, the question remains how prevalent the holonomy representations are. There are clearly none for $\E(\rho)[S] = 0$; this contradicts Gauss-Bonnet. In \cite{Tan94} Tan gives a representation of a genus 3 closed surface $S$ with Euler class $\E(\rho)[S] = 2$, which is not the holonomy of any cone-manifold structure; but he also finds representations arbitrarily close to this one, which do give branched hyperbolic structures.
\begin{qn}
    For a given integer $m \neq 0$, $\chi(S)+1 \leq m \leq -\chi(S)-1$, are
    holonomy representations dense, or conull, in the set of
    representations with $\E(\rho)[S]=m$?
\end{qn}

\subsection{Structure of this paper}

This paper is organised as follows. 

In section \ref{sec:background} we present brief background required for the proofs; this is in addition to the prequel \cite{Me10MScPaper1} and only consists of material that was not required there. We discuss geometric cone-manifold structures and the Euler class of a representation, representation and character varieties, the symplectic structure and measure on the character variety, and the action of the mapping class group and related groups on it.

In section \ref{sec:hyperbolize_pants} we prove Goldman's theorem in the case of pants, which is a building block for the proof in general. This proof, like the results of the prequel, deduces the geometric arrangement of isometries from algebraic data of a representation, and constructs an explicit fundamental domain. Then in section \ref{sec:Goldmans_theorem} we apply our methods to give a proof of Goldman's theorem \ref{goldman_theorem}. Using building blocks of punctured tori and pants, we piece together developing maps to obtain a hyperbolic structure on a larger surface.

In section \ref{sec:constructions_for_genus_2} we turn to the closed genus 2 surface and prove theorem \ref{restricted_version}. We use the Euler class and $\widetilde{PSL_2\R}$ to classify the possible splittings of the surface into two punctured tori. We find two pentagonal fundamental domains which fit
together. 

Finally in section \ref{sec:one_off_extremal} we prove theorem \ref{technical_prop}, as made precise in theorem \ref{thm:technical_precise}. We use Goldman's ergodicity results in \cite{Goldman03}. These allow us to change basis to ``almost go almost anywhere" in the character variety, and hence, simply by changing basis, alter the
geometric situation almost entirely as we please. This is a technique which we hope will have applications to other results.

\subsection{Acknowledgments}

This paper forms one of several papers arising from the author's Masters thesis \cite{Mathews05}, completed at the University of Melbourne under Craig Hodgson, whose advice and suggestions have been highly valuable. It was completed during the author's postdoctoral fellowship at the Universit\'{e} de Nantes, supported by a grant ``Floer Power'' from the ANR.

\section{Background}
\label{sec:background}

Throughout this paper, $S$ denotes a connected oriented surface of finite genus and with finitely many boundary components or cusps.

\subsection{Geometric structures and the Euler class}
\label{sec:euler_class}

Recall that a geometric $(X, \Isom X)$ structure on a surface $S$ can be considered as a metric on $S$ locally isometric to $X$; equivalently, as an atlas of coordinate charts with transition maps; equivalently, as a developing map $\D: \tilde{S} \To X$ equivariant under the action of the fundamental group \cite{Thurston_book, Thurston_notes}. A loop in $S$ gives rise to an isometry in $X$, hence the holonomy representation $\rho: \pi_1(S) \To \Isom X$.

For a given homomorphism $\rho: \pi_1(S) \To \Isom X$, a geometric structure with holonomy $\rho$ can be described as a type of section of a certain fibre bundle: see e.g. \cite{Leleu} for details. Let $\F(S, X, \rho)$ be the flat $X$-bundle over $S$ with holonomy $\rho$, i.e. the quotient of $\tilde{S} \times X$ by $\pi_1(S)$, where $\pi_1(S)$ acts on $\tilde{S}$ by deck transformations, and on $X$ via the isometries given by $\rho$. The product $\tilde{S} \times X$ is foliated by lines of the form $\tilde{S} \times \{x\}$, for each individual $x \in X$; this foliation descends to $\F(S, X, \rho)$. A section $s$ of the bundle transverse to this foliation is precisely a geometric structure on $S$; $s$ immediately gives a developing map by lifting to a map $\tilde{s}: \tilde{S} \To \tilde{S} \times X$, which is equivariant under the action of $\pi_1(S)$, and projecting onto the second coordinate. Transversality of $s$ to the foliation is equivalent to $\tilde{s}$ always having nonzero derivative in the $X$ coordinate, i.e. the developing map being an immersion.

Restrict now to 2-dimensional hyperbolic geometry, i.e. $S$ a surface with $\chi(S) < 0$ and $(X, \Isom X)  = (\hyp^2, PSL_2\R)$. Consider the associated principal bundle $\F(S,PSL_2\R,\rho)$, which is the flat $PSL_2\R$-bundle over $S$ with holonomy $\rho$. The Euler class of $\rho$ arises naturally as an obstruction to finding sections of this bundle. Recall that, since a hyperbolic isometry is determined by the image of a unit tangent vector, $PSL_2\R \cong UT\hyp^2$, the unit tangent bundle of the hyperbolic plane.

Consider a cell complex structure on $S$; we attempt to find a section of $\F(S,PSL_2\R,\rho)$ over $S$, skeleton by skeleton. A section can trivially be found on the 0-skeleton, choosing unit points and tangent vectors in $\hyp^2$ above every vertex, equivariantly, and all such sections are homotopic. We can then extend to a section $s_1$ over the 1-skeleton, joining unit tangent vectors by paths; since $\pi_1(UT\hyp^2) = \Z$, there are infinitely many choices for the extension along each edge; along each edge the unit tangent vectors may spin arbitrarily many times. 

This leaves only extension over the 2-skeleton, which is impossible, since it would give a nowhere zero unit vector field on $S$ and $\chi(S)<0$. However, over a 2-cell $\sigma$ of $S$, $s_1$ gives a unit tangent vector field around $\partial \sigma$, hence a loop in $UT\hyp^2$. The homotopy class of this loop corresponds to the number of times the unit tangent vector ``spins'' as it travels around the loop; we assign this integer to $\sigma$. his gives a 2-cochain of $S$, hence a cocyle; adjustment by a coboundary corresponds to altering the amount of ``spin'' chosen along each particular edge. The cohomology class of this cochain does not depend on the choice of 1-section, nor the cellular decomposition of $S$: it gives the \emph{Euler class} $\E(\rho) \in H^2(S)$ of $\F(S, PSL_2\R, \rho)$, as in \cite{Milnor_Stasheff}.

If $S$ has boundary, we can obtain a \emph{relative Euler class}, which is defined in the same way: $\E(\rho)$ measures the obstruction to extending a section of $\F(S,PSL_2\R,\rho)$ over the 2-skeleton of $S$. It is still independent of choice of triangulation and 1-section; however we first require a trivialization over the boundary. When $S$ has no boundary, different choices of trivializations over edges ``cancel out'' since there are faces on both sides of an edge!

Let $C_i$ denote the boundary curves of $S$. A canonical trivialization over $\partial S$ exists whenever $\rho$ takes each boundary curve to a non-elliptic element of $PSL_2\R$. Since all loops in $\pi_1(S)$ freely homotopic to a boundary curve $C_i$ are conjugate, it makes sense to speak of $\rho(C_i)$ as elliptic, hyperbolic, etc. By abuse of notation, take $C_i \in \pi_1(S)$ some such loop; let $\rho(C_i) = c_i$. A trivialization over the boundary will be given by choosing a preferred lift $\tilde{c_i} \in \widetilde{PSL_2\R}$ of each $c_i \in PSL_2\R$; see \cite{Me10MScPaper0, Me10MScPaper1} for discussion of $\widetilde{PSL_2\R}$. Considered under a developing map, a trivialization over an edge corresponds to choosing a path of unit tangent vectors connecting sections over vertices; around a boundary component of $S$, this corresponds to choosing a path of unit tangent vectors between unit tangent vectors related by the isometry $c_i$; hence to choosing a lift $\tilde{c_i}$ of $c_i$. As discussed in \cite{Goldman88, Me10MScPaper0, Me10MScPaper1}, non-elliptic elements of $PSL_2\R$ have a preferred ``simplest lift'' in $\widetilde{PSL_2\R}$ as the path in $PSL_2\R$ given by restricting to $[0,1]$ the unique homomorphism $\phi: \R \To PSL_2\R$ with $\phi(1) = c_i$; hence we will obtain a canonical trivialization in this case. When each $c_i$ is non-elliptic, then, define $\E(\rho) \in H^2(S)$ to be the relative Euler class with this canonical trivialization.

\begin{figure}
\centering
\includegraphics[scale=0.6]{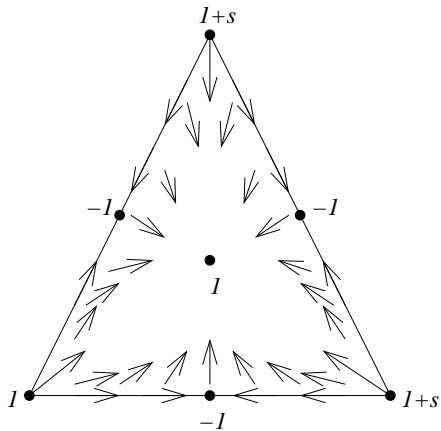}
\caption{The vector field $\V$ on a triangle, with singularities.}
\label{fig:21}
\end{figure}

Now suppose that $\rho$ is the holonomy of a hyperbolic cone-manifold structure on $S$. Suppose there are no corner points, only interior cone points $p_1, \ldots, p_k$, with orders $s_1, \ldots, s_k$; recall as in \cite{Me10MScPaper1}, following \cite{Troyanov}, the \emph{order} $s_i$ of a cone point $p_i$ with cone angle $\theta_i$ is given by $\theta_i = 2\pi(1+s_i)$ for an interior cone point and $\theta = 2\pi( \frac{1}{2} + s_i )$ for a corner point. Take a triangulation of $S$ by geodesic hyperbolic triangles, where all cone points are vertices of the triangulation. There is a ``standard" unit vector field $\V$ on $S$ with precisely one singularity for every vertex, edge and face of $S$. The orders of the singularities are: $1$ on every face; $-1$ on every edge; and $1+s_i$ at every vertex, where $s_i$ is the order of the cone/corner point there (possibly zero). See figure \ref{fig:21}. The sum of the indices of the singularities of $\V$ is then $\chi(S) + \sum s_i$. Thus $\E(\rho)[S] = \pm(\chi(S) + \sum s_i)$, the sign depending on orientation: one way to see this is to choose a different triangulation where singularities occur off the 1-skeleton, so that the spin around every 2-cell is clear; this requires that there are no singularities on the boundary, i.e. corner points.

This argument works, in a limit, even if a boundary component $C_i$ of $S$ is \emph{cusped}, i.e. the section over $C_i$ is a point at infinity. This is possible when $c_i$ is parabolic; the developing image of $C_i$ is $\Fix c_i$; there are no corner points on a cusped boundary.

The above discussion is summarised by the following proposition.
\begin{prop}
\label{euler_class_of_cone_manifold}
    Let $S$ be a surface with boundary. Suppose $\rho: \pi_1(S) \To PSL_2\R$ takes every boundary curve to a non-elliptic element. Suppose $\rho$ is the holonomy of a hyperbolic cone-manifold structure on $S$ with no corner points, each boundary component totally geodesic or cusped, and interior cone points of orders $s_i$. Then the relative Euler class $\E(\rho)$ satisfies
    \[
        \E(\rho)[S] = \pm \left( \chi(S) + \sum s_i \right).
    \]
    \qed
\end{prop}
Since each $s_i > 0$, and since $\sum s_i < -\chi(S)$ (lemma 2.1 of \cite{Me10MScPaper1}), it follows that $|\E(\rho)[S]| \leq |\chi(S)|$. If $S$ has no cone points, then $\E(\rho)[S] = \pm \chi(S)$, as it should, and takes its extremal value. As more cone points are introduced, the Euler class deviates further from this extremal value. 

If we consider several $S_1, \ldots, S_n$ glued together into a larger surface $S$, then we see that the spins along the common boundary cancel out, so that the relative Euler class is additive. For the representation on each piece, we must choose a basepoint $p_i \in S_i$, connected to the basepoint $p$ of the overall surface $S$ by a particular path; then we may define a representation $\rho_i$ on $\pi_1(S_i,p_i)$ and obtain the following result. See also \cite{Goldman_thesis, Goldman84}.
\begin{lem}
\label{euler_additivity}
    Suppose a surface $S$ is decomposed along disjoint simple closed curves $D_i$, with each
    $\rho(D_i)$ not elliptic, into surfaces $S_1, S_2, \ldots, S_n$. Suppose that $\E(\rho)$ is well-defined, i.e. $\rho(C_i)$ is non-elliptic for each boundary component $C_i$ of $S$.
    Then
    \[
        \E(\rho_1)[S_1] + \cdots + \E(\rho_n)[S_n] = \E(\rho)[S].
    \]
    \qed
\end{lem}

\subsection{Algebraic description of the Euler class}
\label{sec:algebraic_description_of_euler_class}

\begin{figure}
\centering
\includegraphics[scale=0.4]{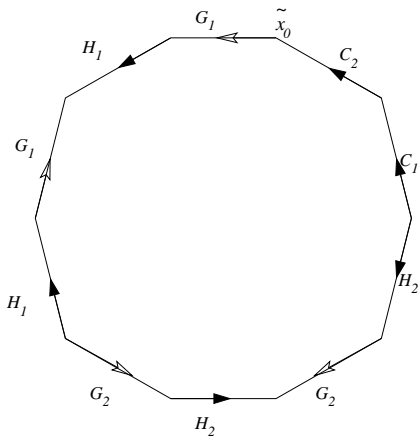}
\caption{Standard cell complex ($k=2$, $n=2$).}
\label{fig:13}
\end{figure}

There is a directly algebraic way to see $\E(\rho)$. Consider the surface $S$ of genus $k$, with $n$ boundary components, and assume that each $c_i$ is non-elliptic. This surface is homotopy equivalent to a standard cell complex $S_0$ with one 0-cell, $2k+n$ 1-cells, and one 2-cell, a $(4k+n)$-gon, glued as shown in figure \ref{fig:13}. We have the standard presentation of the fundamental group $\pi_1(S)$; the relator says that this $(4k+n)$-gon bounds a disc.
\[
    \Bigl\langle G_1, H_1, \ldots, G_k, H_k, C_1,
    \ldots, C_n \; | \; [G_1, H_1] \cdots [G_k, H_k] C_1 C_2
    \cdots C_n = 1 \Bigr\rangle.
\]
Let $\rho(G_i) = g_i$, $\rho(H_i) = h_i$, $\rho(C_i)=c_i$; we arbitrarily choose lifts $\tilde{g}_i, \tilde{h}_i \in \widetilde{PSL_2\R}$; since all the $c_i$'s are non-elliptic, we have a preferred lift $\tilde{c}_i$ of each. 

A partial section $s_1$ of $\F(S_0, PSL_2\R, \rho)$ over the 1-skeleton of $S_0$ gives a loop in $UT\hyp^2$ over the boundary of our polygon. The paths of unit tangent vectors over the edges of the polygon give elements of $\widetilde{PSL_2\R}$; we can take a section such that the paths along various edges are given by (appropriate conjugates of) the $\tilde{g}_i, \tilde{h}_i, \tilde{c}_i$. Moving anticlockwise around the polygon in $\tilde{S}$, we thus obtain a loop in $UT\hyp^2$ which is represented by
\[
    [\tilde{g}_1, \tilde{h}_1] \cdots [\tilde{g}_k, \tilde{h}_k] \;
    \tilde{c}_1 \cdots \tilde{c}_n.
\]
This is a path of vectors which spins some number $m$ of times, i.e. is equal to $\z^m$; it is a lift of $[g_1, h_1] \cdots [g_k, h_k] c_1 \cdots c_n = 1$. Recall (e.g. lemma 2.6 of \cite{Me10MScPaper1}) that $[g_i, h_i]$ lifts to $\widetilde{PSL_2\R}$ independent of choices of lifts $\tilde{g}_i, \tilde{h}_i$. By the discussion of the Euler class above, we immediately have the following.

\begin{prop}
\label{euler_algebra}
    Let $S$ be an orientable surface with $\chi(S) < 0$. Let $\rho: \pi_1(S) \To PSL_2\R$ be a representation, and let $\pi_1(S)$ have the presentation given above, where no $c_i$ is elliptic. The (possibly relative) Euler class $\E(\rho)$ takes the fundamental class $[S]$ to $m \in \Z$ where the unique lift of the relator
    \[
    [\tilde{g}_1, \tilde{h}_1] \cdots [\tilde{g}_k,
    \tilde{h}_k] \; \tilde{c}_1 \cdots \tilde{c}_n,
    = \z^m.
    \]
    \qed
\end{prop}
Thus the Euler class is an obstruction to lifting $\rho$ to $\widetilde{PSL_2\R}$. The \emph{Milnor-Wood inequality} \cite{Milnor, Wood, Me10MScPaper1} states that for such a product, $|m| \leq |\chi(S)|$. For holonomy representations, this follows immediately from the above and proposition \ref{euler_class_of_cone_manifold}. A proof using our methods of twisting, is given in our \cite{Me10MScPaper0}; see that paper for general comments on the inequality.

When $\rho$ is \emph{abelian}, then all of the $g_i,h_i,c_i$ are all hyperbolics with the same axis, elliptics with the same fixed point, or parabolics with the same fixed point at infinity. Thus we easily see the relator is $1 \in \widetilde{PSL_2\R}$ and hence $\E(\rho) = 0$.

When $S$ is a \emph{punctured torus}, writing $g=g_1$, $h=h_1$, $c=c_1$, $\E(\rho)$ is well-defined iff $[g,h]$ is non-elliptic. Recall the theorem in the prequel (2.9 of \cite{Me10MScPaper1} or 3.14 of \cite{Me10MScPaper0}) giving the possible regions of $\widetilde{PSL_2\R}$ in which a commutator may lie (also appearing in \cite{Wood, Eisenbud-Hirsch-Neumann, Goldman88}). We use this to prove the following.
\begin{prop}
\label{punctured_torus_euler_class}
    Let $S$ is a punctured torus and assume the relative Euler class is
    well-defined. Then:
    \begin{enumerate}
        \item
            $\Tr[g,h] \geq 2$ is equivalent to $\E(\rho)[S] = 0$;
        \item
            $\Tr[g,h] \leq -2$ is equivalent to $\E(\rho)[S] = \pm 1$. Moreover in this case $\E(\rho)[S] = +1$ or $- 1$ as $[g,h] \in \Hyp_1 \cup \Par_1$ or $\Hyp_{-1} \cup \Par_{-1}$ respectively.
    \end{enumerate}
\end{prop}

\begin{Proof}
    As $\E(\rho)$ is well-defined, $[g,h]$ is not elliptic, and $\Tr[g,h] \in (-\infty, -2] \cup [2, \infty)$.

    Suppose $\Tr[g,h] \geq 2$. Then $[g,h]$ in $\{1\}$, $\Hyp_0$ or $\Par_0$. Now $c^{-1} = [g,h]$, hence $\tilde{c}$, the simplest lift of $c$, satisfies $\tilde{c}^{-1} = [g,h]$, so $[g,h] \tilde{c} = 1 \in \widetilde{PSL_2\R}$. Hence $\E(\rho)[S] = 0$ by proposition \ref{euler_algebra} above.

    Now suppose $\Tr[g,h] \leq -2$. Then $[g,h]$ lies in $\Hyp_{\pm 1}$, $\Par_{-1}^+$ or $\Par_1^-$. Assume $[g,h] \in \Hyp_1 \cup \Par_1^-$; the reversed-orientation case is similar. Since $c^{-1} = [g,h]$, then the simplest lift $\tilde{c}$ lies in $\Hyp_0 \cup \Par_0^-$ and satisfies $\tilde{c}^{-1} = \z^{-1} [g,h]$. Thus $[g,h] \tilde{c} = \z$ and $\E(\rho)[S] = 1$.
\end{Proof}

\subsection{The character variety, its symplectic structure and measure}
\label{sec:character_variety}

In section 4 of \cite{Me10MScPaper1} we gave some discussion of representation and character varieties for a punctured torus. Here we extend that discussion to more general surfaces. Recall the representation variety $R_G(S)$ describes all homomorphisms $\rho: \pi_1(S) \To G$; for our purposes $G = SL_2\R$ or $PSL_2\R$ (although much of the following is true for much more general $G$, see \cite{Goldman84}); $R_G(S)$ is a closed algebraic set. For a closed surface $S$ of genus $g \geq 2$, $R_{PSL_2\R}(S)$ is not connected. In fact if we vary a representation continuously, $\E(\rho)[S]$ changes continuously, but is an integer, hence remains constant. For \emph{closed} surfaces, Goldman classified the components of $R(S)$ completely.
\begin{thm}[Goldman \cite{Goldman88}]
\label{topological_components}
    For a closed surface $S$ with $\chi(S) < 0$, $R_{PSL_2\R}(S)$ has precisely $2|\chi(S)|+1$
    components, parameterised by the Euler class.
    \qed
\end{thm}

Although every $SL_2\R$ representation projects to a $PSL_2\R$ representation, not every $PSL_2\R$ representation lifts to an $SL_2\R$ representation. Taking a standard presentation for $\pi_1(S)$ with one relator, in a $PSL_2\R$ representation the relator gives a product in $PSL_2\R$ multiplying to the identity. Choosing lifts of images of generators to $SL_2\R$, this relator multiplies to $1$ or $-1$; when $S$ is a \emph{closed} surface, the relator is a product of commutators, which lifts \emph{uniquely} to $SL_2\R$ as $\pm 1$. However we may consider \emph{twisted representations} into $SL_2\R$ as in \cite{Goldman84, Huebschmann}, allowing the relator to multiply to $\pm 1 \in SL_2\R$. For any surface $S$ then we obtain the twisted representation spaces $R_{SL_2\R}^{\pm}(S)$, where $R_{SL_2\R}^{+}(S) = R_{SL_2\R}(S)$. Each $PSL_2\R$ representation lifts to a twisted $SL_2\R$ representation, so $R_{PSL_2\R}(S)$ is an obvious quotient of $R_{SL_2\R}^{\pm}(X) = R_{SL_2\R}(S) \cup R_{SL_2\R}^{-} (S)$.

For a representation $\rho: \pi_1(S) \To SL_2\R$, the character $\chi$ of $\rho$ is the function $\Tr \circ \rho: \pi_1(S) \To \R$. It is determined by its values at finitely many elements $\gamma_1, \ldots, \gamma_m \in \pi_1(S)$; taking $t:R_{SL_2\R}(S) \To \R^m$, $\rho \mapsto (\Tr \circ \rho (\gamma_1), \ldots, \Tr \circ \rho (\gamma_m) )$, the character variety is $X(S) = t(R_{SL_2\R}(S)) \subset \R^m$; it is again a closed algebraic set \cite{Culler-Shalen}. Recall $SL_2\R$ acts on the representation space by conjugation, and the quotient by this action can be identified with $X(S)$ away from singularities. For twisted $SL_2\R$ representations we may do the same, taking the trace of the image of the same generators, and obtain the character variety $X^-(S)$ of representations in $R_{SL_2\R}^{-}(S)$ respectively; this may be identified with $R_{SL_2\R}^-(S) / SL_2\R$ away from singularities. Since both $X(S)$ and $X^-(S)$ are obtained by taking traces of images of various generators, they may be regarded as disjoint subsets of the same $\R^m$. We write $X^\pm(S) = X(S) \cup X^-(S)$. For $PSL_2\R$ represenations, we may lift to a \emph{possibly twisted} $SL_2\R$ representation and again take traces to obtain a character. Thus we may speak of a character in $X^{\pm}(S)$ of a $PSL_2\R$ representation.

For a general surface with boundary, and $G = SL_2\R$ or $PSL_2\R$ (or any of a large class of Lie groups \cite{Goldman84}), there is a symplectic structure on $R_G(S)/G$, although the structure is singular along the singularities of $R_G(S)$ (\cite{Goldman97}): the same is also true for twisted representations \cite{Huebschmann}. We briefly describe how this structure arises: see \cite{Hodgson_thesis} or \cite{Goldman84} for more details. Consider a smooth path $\rho_t$ of representations in $R_G(S)$. Approximate $\rho_t$ to first order:
\[
    \rho_t (x) = \exp \left( t u(x) + O(t^2) \right) \rho_0(G),
\]
where $u$ is some function $\pi_1(S) \To \mathfrak{g}$, the Lie algebra of $G$. Since each $\rho_t$ is a homomorphism,
\[
    u(xy) = u(x) + \Ad (\rho_0(x)) \; (u(y)),
\]
where $\Ad: G \To \Aut (\mathfrak{g})$ is the adjoint representation. There is in fact an $\R\pi_1(S)$-module structure on $\mathfrak{g}$, given by $x.v = \Ad(\rho_0 (x)) (v)$, for $x \in \pi_1(S)$ and $v \in \mathfrak{g}$. We denote this $\R\pi_1(S)$-module $\mathfrak{g}_{\Ad \rho_0}$. The condition for $u$ above is then just $u(xy) = u(x) + x.u(y)$, i.e. that $u: \pi_1(S) \To \mathfrak{g}_{\Ad \rho_0}$ is a \emph{1-cocycle} in the group cohomology of $\pi_1(S)$ with coefficients in $\mathfrak{g}_{\Ad \rho_0}$ (see e.g. \cite{Brown} or \cite{HS} for details). The Zariski tangent space to $R_G(S)$ at $\rho_0$ can be identified with the $\R$-vector space structure on these cocycles.

Now consider the tangent space of $R_G(S)/G$, the quotient space by conjugation. A path $\rho_t \in R_G(S)$ of representations corresponds to a constant path in this quotient if and only if, to first order, each $\rho_t$ is conjugate to $\rho_0$, i.e. $\rho_t (x)= g_t^{-1} \rho_0(x) g_t$ for some path $g_t \in G$. So let $g_t = \exp(tu_0 + O(t^2))$ and again let $\rho_t = \exp \left( t u(x) + O(t^2) \right) \rho_0$. The condition that $u$ give a constant path is precisely the coboundary condition
\[
    u(x) = \Ad (\rho_0(x)) (u_0) - u_0 = \delta u_0 \in B^1 \left( \pi_1(S); \mathfrak{g}_{\Ad \rho_0} \right).
\]
So the tangent space to $R_G(S)/G$ at $[\rho_0]$ is the (vector-space structure on the) cohomology module $H^1(\pi_1(S); \mathfrak{g}_{\Ad \rho_0})$. Note that this group cohomology module is also $H^1(S;\B)$, for a bundle $\B$ of coefficients over $S$ associated with the $\pi_1(S)$-module $\mathfrak{sl}_2\R_{\Ad \rho}$.

For \emph{closed} surfaces, the dimensions of these tangent spaces are given in \cite{Goldman84}. The dimension of the tangent space to $R_G(S)$ at $\rho_0$ is $\dim Z^1(\pi_1(S);  \mathfrak{sl}_2\R_{\Ad \rho_0}) = 6g-3+ \dim C(\rho_0)$, where $C(\rho_0)$ is the centraliser of $\rho_0(\pi_1(S))$ in $G$. And $\dim B^1 (\pi_1(S); \mathfrak{sl}_2\R_{\Ad \rho_0}) = 3 - \dim C(\rho_0)$. Thus the dimension of the tangent space to $R_G(S)/G$ at $[\rho_0]$ is $\dim H^1(\pi_1(S); \mathfrak{sl}_2\R_{\Ad \rho_0}) = 6g-6 + 2 C(\rho_0)$. This $C(\rho_0)$ is trivial for non-abelian $\rho_0$, $1$-dimensional for non-trivial abelian $\rho_0$, and all of $G$ (hence $3$-dimensional) for $\rho_0=1$. Letting $R_G(S)^-$ denote the non-abelian representations, we may take the quotient $R_G(S)^-/G$, which is
$(6g-6)$-dimensional. In general however this space is not Hausdorff: \cite{Goldman84}. The characters of abelian representations are precisely the singularities of $R_G(S)/G$.

Returning to general surfaces, consider the cup product in group cohomology on $\pi_1(S)$ with coefficients in $\mathfrak{g}_{\Ad \rho_0}$. This gives a dual pairing
\[
    H^1 \left( \pi_1(S); \mathfrak{g}_{\Ad \rho_0} \right) \times
    H^1 \left( \pi_1(S); \mathfrak{g}_{\Ad \rho_0}^* \right) \To
    H^2(\pi_1(S); \R) = \R.
\]
Since the Killing form on $\mathfrak{psl}_2\R \cong \mathfrak{sl}_2\R$ is nondegenerate and invariant under the adjoint representation, there is an isomorphism $\mathfrak{g}_{\Ad \rho_0} \cong \mathfrak{g}_{\Ad \rho_0}^*$. Using this isomorphism with the cup product we can define a dual pairing on $R_G(S)/G$
\[
    \omega_{\rho_0} : \; H^1 \left( \pi; \mathfrak{g}_{\Ad \rho_0} \right) \times H^1 \left( \pi; \mathfrak{g}_{\Ad \rho_0} \right) \To \R.
\]
From the above, $\omega_{\rho_0}$ is actually a 2-form on $T_{[\rho_0]} \left( R_G(S)/G \right)$ (i.e. on the tangent space at $\rho_0$). This clearly varies continuously with $[\rho_0]$, so we obtain a 2-form $\omega$ on $R_G(S)/G$, which is singular at the equivalence classes of abelian representations. It can be shown (see \cite{Goldman84}) that $\omega$ is closed and nondegenerate.

If $S$ is a closed surface, then $R_G(S)/G$ is everywhere even-dimensional (even though the dimension varies) and we obtain a symplectic structure on $R_G(S)/G$. Hence we obtain a symplectic structure on $X(S)$, away from the characters of abelian representations. By taking an appropriate exterior power of $\omega$, we obtain an area form on $R_G(S)^-/G$, and a singular area form on $R_G(S)/G$. This gives a measure on $R_G(S)/G$. It can be shown that the singular set has measure zero and the measure of $R_G(S)/G$ is finite: see \cite{Goldman97, Huebschmann}. This is also true for twisted representations \cite{Huebschmann}. So we obtain a measure $\mu_S$ on $X^{\pm}(S)$. Considering $X^{\pm}(S)$ as a subset of some $\R^{m}$, away from singular points the top power of $\omega$ is some multiple of the Euclidean area form, hence $\mu_S$ is absolutely continuous with respect to Lebesgue measure.

Having defined $\mu_S$ and explained the details of characters of $PSL_2\R$ representations, we can refine theorem \ref{technical_prop} to a precise statement.
\begin{thm}
\label{thm:technical_precise}
    Let $S$ be a closed orientable surface of genus $g \geq 2$. Let $\Omega_S \subset X^{\pm}(S)$ denote the set of characters of representations $\rho: \pi_1(S) \To PSL_2\R$ such that
    \begin{enumerate}
        \item
            $\E(\rho)[S] = \pm (\chi(S) + 1)$;
        \item
            there exists a non-separating simple closed curve $C$ on $S$ such
            that $\rho(C)$ is elliptic.
    \end{enumerate}
    Then $\mu_S$-almost every character in $\Omega_S$ is the character of a holonomy representation for a cone-manifold structure on $S$ with a single cone point with
    cone angle $4\pi$.
\end{thm}
(In fact, since $\E(\rho)[S]$ is odd, $\rho$ only lifts to twisted representations in $R_{SL_2\R}^{-}(S)$ and hence has character in $X^{-}(S)$.)

Although the above considers surfaces, the character variety can be defined in a similar way for any manifold. For a circle $S^1$ we obtain $X(S^1) \cong \R$. At points other than $\pm 2$ the character defines the conjugacy class of a representation uniquely.

For a surface $S$ with boundary, we can consider a \emph{relative character variety}, following \cite{Goldman97}. The boundary $\partial S$ is a collection of circles $C_1, \ldots, C_n$, and so $X(\partial S) = X(S^1)^n = \R^n$. There is then a restriction map
\[
    \partial^\#: X(S) \To X(\partial S) = \R^n.
\]
If we specify for each $C_i$ a conjugacy class $\mathcal{C}_i$, then we may define the \emph{relative character variety} to be
\[
    X_\mathcal{C}^{\pm} (S) = \left\{ [\rho] \in R_{SL_2\R}^{\pm}(S)/SL_2\R \; | \; \rho(C_i) \in \mathcal{C}_i \right\}.
\]
Note that if $\mathcal{C}_i$ is hyperbolic or elliptic, then it is described completely by its trace $t$ (while a trace of $\pm 2$ is ambiguous), and we can write $X_t(S)$. This agrees with our notation for the relative character variety of the punctured torus in the prequel \cite{Me10MScPaper1}.

Starting from a \emph{closed} surface $S$, if we \emph{cut} $S$ along a curve $\alpha$ to obtain a surface $(S|\alpha)$ then the inclusion $(S|\alpha) \to S$ induces a map $X^\pm(S) \to X^\pm(S|\alpha)$. Letting $s_\alpha: R_{SL_2\R}^\pm(S) / SL_2\R \To \R$ be defined by $s_\alpha [\rho] = \Tr(\rho(\alpha))$, we may disintegrate the measure $\mu_S$ on $X^\pm(S)$ to obtain a measure $\mu_t$ on each $s_\alpha^{-1}(t) \subset X_t(S)$. See \cite{Goldman97} for further details.

When $S$ is a punctured torus, we have described the character variety $X(S)$ and relative character variety $X_t(S)$ in section 4.2 of \cite{Me10MScPaper1}, following \cite{Goldman88}. By theorem 4.1 there, $X(S)$ is the set of all $(x,y,z)$ with $\kappa(x,y,z) \geq 2$ or at least one of $|x|, |y|, |z| \geq 2$; here $(x,y,z) = (\Tr g, \Tr h, \Tr gh)$ and $\kappa(x,y,z) = \Tr[g,h] = x^2 + y^2 + z^2 - xyz - 2$. Then $X_t(S) = X(S) \cap \kappa^{-1}(t)$, which is $2$-dimensional; we obtain a measure, indeed a symplectic structure, on each $X_t(S)$ and the symplectic form can be written explicitly: see \cite{Goldman03}.

\subsection{The action on the character variety}
\label{sec:action_on_character_variety}

We now consider the effect of changing a (possibly twisted) representation $\rho: \pi_1(S) \To SL_2\R$ by pre-composition with an automorphism of $\pi_1(S)$: that is, take $\phi \in \Aut \pi_1(S)$ and replace $\rho$ with $\rho' = \rho \circ \phi$. This descends to an action of $\Aut \pi_1 (S)$ on the character variety. Since traces are invariant under conjugation, the action of $\Inn \pi_1(S)$ is trivial and we consider the action of the quotient $\Out \pi_1(S) = \Aut\pi_1(S)/\Inn \pi_1(S)$. Points in $X^\pm (S)$ which are related under this action ought to be considered as equivalent in terms of the underlying geometry.

In the prequel \cite{Me10MScPaper1} we described $X(S)$, and the orbits of $X(S)$ under this action, precisely, for $S$ a punctured torus. Recall (proposition 4.6) that characters of irreducible representations $(x,y,z)$, $(x',y',z') \in X(S)$ are equivalent iff they are related by some sequence of the moves $(x,y,z) \mapsto (x,y,xy-z)$, $(x,y,z) \mapsto (-x,-y,z)$ and permutations of coordinates. These are called \emph{Markoff} triples. The equivalence relation can be considered as the action of a semidirect product $\Gamma = PGL_2\Z \ltimes \left( \frac{\Z}{2} \oplus \frac{\Z}{2} \right)$ on $X(S)$. This action preserves each relative character variety $X_t(S)$. For more detail see \cite{Goldman03}.

Also recall the geometric interpretation $\Out \pi_1(S) \cong \mcg(S)$, when $S$ is a closed surface or a punctured torus, the Dehn--Nielsen theorem \cite{Stillwell, Nielsen27, Goldman03}.

For $S$ closed, the 2-form $\omega$ is invariant under the action of $\Out \pi_1 (S) \cong \mcg (S)$ on $X^\pm(S)$, and hence the action of $\Out \pi_1(S)$ is measure-preserving with respect to $\mu_S$ \cite{Goldman97}. For $S$ with boundary, the action preserves the measure on each relative character variety $X_t(S)$. In particular this is true for the punctured torus, where $X_t(S) = X(S) \cap \kappa^{-1}(t)$.

\section{How to hyperbolize your pants}
\label{sec:hyperbolize_pants}

In the prequel \cite{Me10MScPaper1} we described how to hyperbolize punctured tori. In order to prove results for higher genus surfaces, we will need to cut them into punctured tori and pants. Thus, while following \cite{Me10MScPaper1} we are masters of punctured tori, we are not yet masters of our pants. No cone points will be considered in this section; we will only need complete hyperbolic structures with totally geodesic or cusped boundary. In this section we prove the following proposition.
\begin{prop}
\label{pants_construction}
    Let $S$ be a pair of pants with $C_1, C_2, C_3 \in \pi_1(S)$ representing boundary curves and $C_1 C_2 C_3 = 1$. Let $\rho: \pi_1(S) \To PSL_2\R$ be a representation taking each boundary curve $C_i$ to a non-elliptic element. Suppose $\E(\rho)[S_i] = 1$ (resp. $-1$). Then $\rho$ is the holonomy of a complete hyperbolic structure on $S$. Each $C_i$ has hyperbolic or parabolic holonomy, and accordingly $C_i$ is totally geodesic or cusped. Each $C_i$ bounds $S$ to its right (resp. left).
\end{prop}

So $\pi_1(S_i) = \langle C_1, C_2, C_3 \; | \; C_1 C_2 C_3 = 1 \rangle$ and let $\rho_i(C_j) = c_j$, each $c_j$ non-elliptic, with simplest lift $\tilde{c}_j$. Since $\E(\rho_i)[S_i] = 1$, by proposition \ref{euler_algebra} above, $\tilde{c}_1 \tilde{c}_2 \tilde{c}_3 = \z$. If some $c_i$ were the identity, so would be $\tilde{c}_i$; hence the other two $\tilde{c}_j$ would be inverses, and could not multiply to $\z$; thus each $c_i$ is hyperbolic or parabolic.

We will use Milnor's angle function $\Theta$ and the twist function to deduce properties of the $\tilde{c}_i$: see \cite{Me10MScPaper0} for details. Alternatively, we could just use the twist function, since $\Theta (\tilde{\alpha}) = 2 \Twist (\tilde{\alpha}, i)$.

First, $\Theta(\tilde{c}_1 \tilde{c}_2 \tilde{c}_3) = \Theta(\z) = \pi$. Thus
\[
  \Theta \left( \tilde{c}_1 \tilde{c}_2 \right) = \Theta \left( \z \tilde{c}_3^{-1} \right) = \pi + \Theta \left( \tilde{c}_3^{-1} \right) = \pi - \Theta \left( \tilde{c}_3 \right).
\]
But since $\tilde{c}_3$ is a simplest lift, $|\Theta(\tilde{c}_3)| \leq \pi/2$, so $\Theta(\tilde{c}_1 \tilde{c}_2) \in (\pi/2, 3\pi/2)$. As $c_1 c_2 = c_3^{-1}$, which is hyperbolic or parabolic, then $\tilde{c}_1 \tilde{c}_2 \in \Hyp_1 \cup \Par_1$: we know $\Theta$ of various regions of $\widetilde{PSL_2\R}$. Thus $\Tr(\tilde{c}_1), \Tr(\tilde{c}_2) \geq 2$ and $\Tr(\tilde{c}_1 \tilde{c}_2) \leq -2$. Hence
\[
  \Tr \left( c_1 \right) \Tr \left( c_2 \right) \Tr \left( c_1 c_2 \right) \leq -8;
\]
this makes sense since lifting the $c_i$ to $SL_2\R$ with either sign does not change the product of traces.

\begin{lem}
    $\Tr[c_1,c_2] > 2$.
\end{lem}

\begin{Proof}

\begin{figure}
\centering
\includegraphics[scale=0.3]{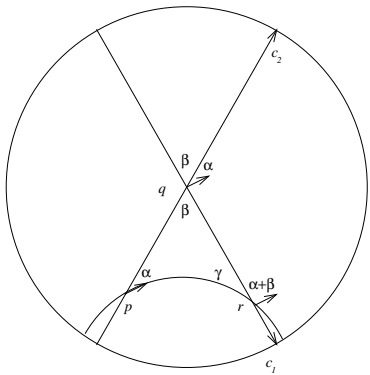}
\caption{Unit vector chase if $\Tr[c_1,c_2] < 2$.} \label{fig:88b}
\end{figure}

Suppose $\Tr[c_1,c_2] < 2$: this is equivalent to $c_1, c_2$ being both hyperbolic with axes intersecting at a point $q \in \hyp^2$ (\cite{Goldman03} or proposition 3.16 of \cite{Me10MScPaper0}). Let $p = c_2^{-1}(q)$ and let $r = c_1(q)$. Let the angles in triangle $pqr$ be $\alpha, \beta, \gamma$ as shown in figure \ref{fig:88b}, so $\alpha + \beta + \gamma < \pi$. Chase unit vectors commencing from the vector at $p$ pointing towards $r$. Under $\tilde{c}_2$ and $\tilde{c}_1$, we obtain the vectors shown, so (taking into account the two possible orientations) $\Twist(\tilde{c}_1 \tilde{c}_2,p) = \pm (\pi - \alpha - \beta - \gamma) \in (-\pi, \pi)$, contradicting $\tilde{c}_1 \tilde{c}_2 \in \Hyp_1 \cup \Par_1$.

If $\Tr[c_1, c_2] = 2$, by proposition 4.2 of the prequel \cite{Me10MScPaper1} (see also \cite{Culler-Shalen, Goldman03}), $c_1, c_2$ form a reducible representation. As $\rho$ is reducible and non-abelian, lemma 4.11 of \cite{Me10MScPaper1} describes $c_1, c_2$: either one of $c_1, c_2$ is hyperbolic and the other parabolic, with a common fixed point; or both $c_1, c_2$ are hyperbolic, with exactly one shared fixed point. In both these cases, a similar unit vector chase contradicts $\tilde{c}_1 \tilde{c}_2 \in \Hyp_1 \cup \Par_1$.
\end{Proof}

\begin{figure}[tbh]
\centering
\includegraphics[scale=0.3]{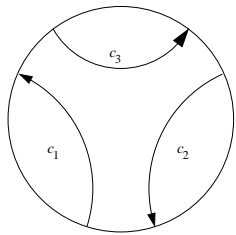}
\caption{Arrangement of axes for an Euler-class-1 pants
representation.} \label{fig:89}
\end{figure}

\begin{Proof}[Of proposition \ref{pants_construction}]
As $\Tr[c_1, c_2]>2$ and $\Tr(c_1) \Tr(c_2) \Tr(c_1 c_2) \leq -8$, we may apply lemma 5.13 of \cite{Me10MScPaper1}. If the $c_j$ are hyperbolic, this lemma tells us that the axes of $c_1, c_2, c_3$ are disjoint and bound a common region, as in figure \ref{fig:89}. If a $c_j$ is parabolic $c_j$, we may take a limit and consider the ``axis" of $c_j$ to degenerate to a point at infinity.

If $c_j$ is hyperbolic then it is the composition of two reflections in lines perpendicular to $\Axis c_j$. If $c_j$ is parabolic then it is the composition of two reflections in lines through its fixed point at infinity. We may take one of these lines to be the common perpendicular of $\Axis c_1$ and $\Axis c_2$, or if $c_j$ is parabolic then we take this line to run to the fixed point at infinity of $c_j$. Then $c_1 c_2 = c_3^{-1}$ is the composition of two reflections. We can do the same for $c_2 c_1$. The axes are as shown in figure \ref{fig:90}.

\begin{figure}[tbh]
\centering
\includegraphics[scale=0.4]{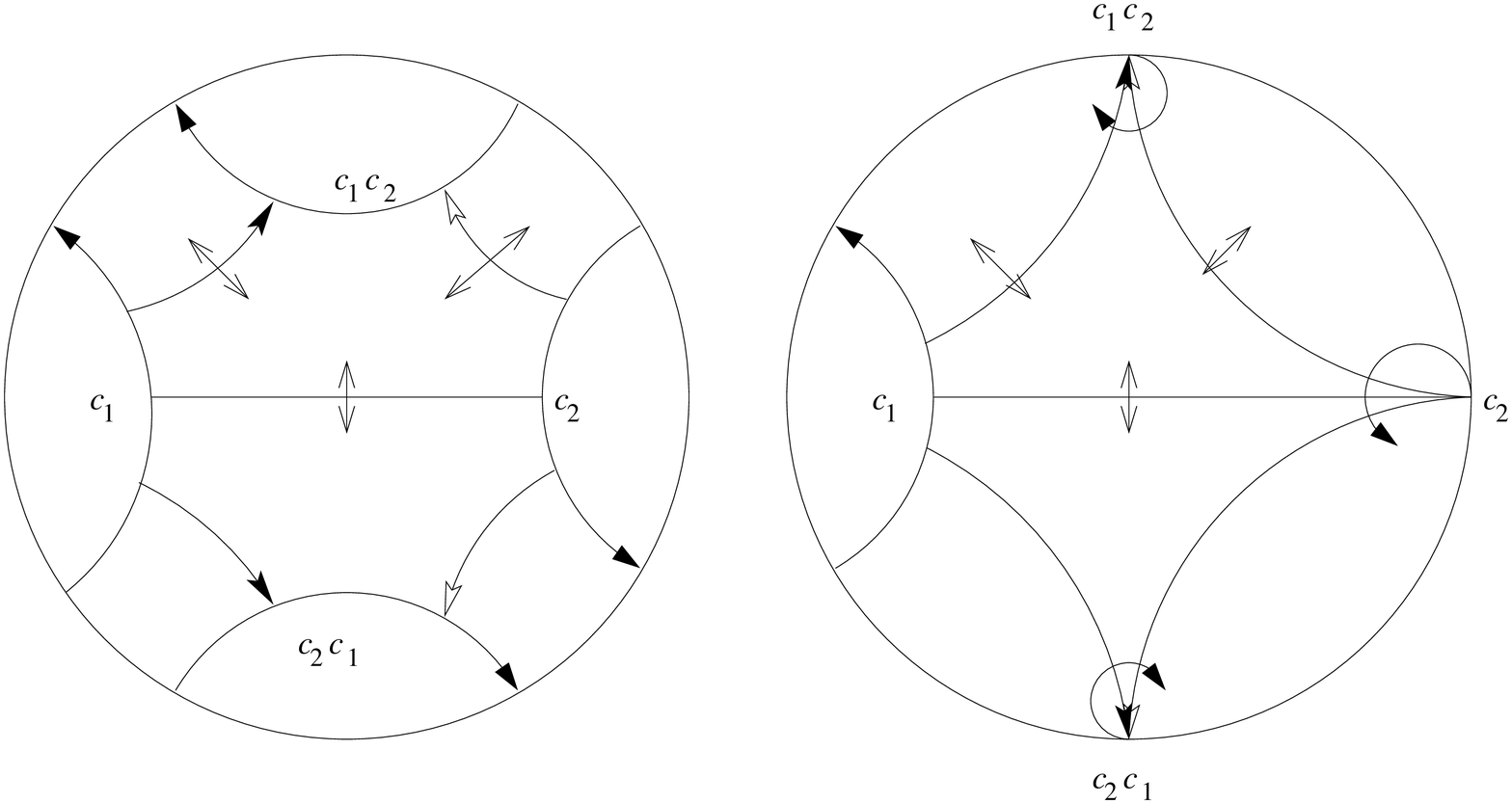}
\caption{Fabricating pants.} \label{fig:90}
\end{figure}

Note that the (possibly degenerate) octagon shown has two pairs of sides identified under $c_1, c_2$, so it forms a fundamental domain for a pair of pants. Since all the angles in the octagon are right angles or $0$, the boundary edges wrap up to give geodesic boundary, or cusps. A developing map is easily then constructed to give a complete hyperbolic structure on $S$ with totally geodesic or cusped boundary, accordingly as each $c_j$ is hyperbolic or parabolic, with holonomy $\rho$. (In fact $\rho$ is discrete and $S$ is the quotient of the convex core of $\rho$.) We also see $S$ is oriented as claimed.
\end{Proof}

\section{Goldman's theorem}
\label{sec:Goldmans_theorem}

In this section we prove Goldman's theorem. Let $S$ be a connected orientable surface of genus $g$ with $n$ boundary components. Assume $\rho$ is non-elliptic on each boundary component, so $\E(\rho)$ is well-defined. Assume $\E(\rho)[S] = - \chi(S) > 0$; the case $\E(\rho)[S] = \chi(S) < 0$ is similar with reversed orientation.

\subsection{Splitting up is hard to do}
\label{sec:splitting}

In \cite{Gallo_Kapovich_Marden}, Gallo--Kapovich--Marden show that for any non-elementary representation $\rho: \pi_1(S) \To PSL_2\C$, where $S$ is a closed oriented surface with $\chi(S)<0$, there exists a system of disjoint curves $C_i$ decomposing $S$ into pants $P_j$, such that the restriction of $\rho$ to each $P_j$ is a 2-generator classical Schottky group. The proof is long, but their methods apply immediately to the case of representations $\pi_1(S) \To PSL_2\R$, and when there are boundary components. That the restriction of $\rho$ to each $P_j$ is a 2-generator classical Schottky group implies that each $\rho(C_i)$ is hyperbolic. An elementary representation can be represented by diagonal matrices, hence has Euler class zero.

The proof of the theorem applies Dehn twists to obtain sufficiently ``complicated" curves that they have holonomy with large trace. Algorithmically, it cuts $g-1$ ``handles", one at a time, so that the genus decreases by $1$ at each stage; and from the remaining piece of genus $1$, cuts off pants (choosing pairs of boundary circles to form into pants arbitrarily each time) until the genus $1$ piece is just a once-punctured torus; then this too is cut into pants. But since, following \cite{Me10MScPaper1}, we are comfortable with punctured tori, we could perform the algorithm so $g$ of the pants have pairs of boundary curves identified, and we glue them back together to give punctured tori. So we can decompose $S$ along curves with hyperbolic holonomy into $g$ tori and $g+n-2$ pants.

\begin{figure}
\centering
\includegraphics[scale=0.35]{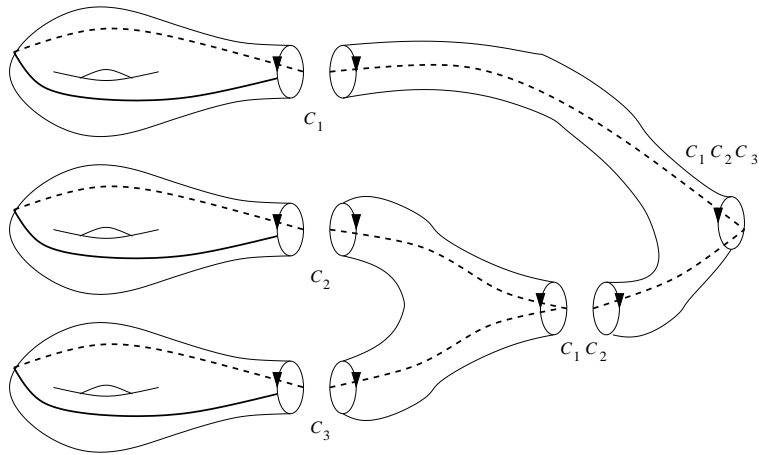}
\caption{Connecting up surfaces, and fundamental group.}
\label{fig:91}
\end{figure}

Then we can assume the surfaces combinatorially fit together as in figure \ref{fig:91}. If $S$ is closed of genus $2$, then we just have two punctured tori. Otherwise, none of the punctured tori are adjacent. We draw all the punctured tori leftmost; these must then be connected together. If $S$ has no boundary, we simply connect up all the punctured tori by pants. If $S$ has boundary, we may add on further pants to the situation of figure \ref{fig:91} to obtain more boundary components.

In short: their theorem implies the following.
\begin{thm}[Gallo--Kapovich--Marden \cite{Gallo_Kapovich_Marden}]
\label{decomposition_thm}
    Let $S$ be an oriented surface with $\chi(S)<0$ and let $\rho: \pi_1(S) \To PSL_2\R$ be a representation with $\E(\rho)[S]$ well-defined and equal to $\pm \chi(S)$. Then there exists a system of disjoint curves $C_i$ decomposing $S$ into pants and punctured tori, such that each $\rho(C_i)$ is hyperbolic.
    \qed
\end{thm}

Let $S_i$ denote the subsurfaces into which $S$ is decomposed. Consider their fundamental groups. We specify basepoints $q_i \in S_i$ and $q \in S$. On each punctured torus we take $q_i$ on the boundary. On each pants we arbitrarily specify a basepoint. We arbitrarily choose one of the $q_i$ to be $q$. To specify how $\pi_1(S_i, q_i)$ and $\pi_1(S,q)$ relate, take a combinatorial tree $T$ dual to the decomposition of $S$, with one vertex for each $q_i$, and a map $\T: T \To S$ mapping the vertices of $T$ onto the corresponding $q_i$. This gives well-defined paths between the $q_i$. We have inclusions $\iota_i: \pi_1(S_i, q_i) \hookrightarrow \pi_1(S, q_i)$ (note basepoints). Let $\alpha_i$ be the unique path from $q$ to each
$q_i$ along the tree $\T$, then we have isomorphisms
\[
    \zeta_i: \pi_1 (S, q_i) \stackrel{\cong}{\To} \pi_1 (S, q),
    \quad x \mapsto \alpha_i.x.\alpha_i^{-1}.
\]
For each $S_i$ we have a representation $\rho_i: \pi_1(S_i,
q_i) \To PSL_2\R$ given by the composition
\[
    \pi_1(S_i, q_i) \stackrel{\iota_i}{\To} \pi_1(S, q_i)
    \stackrel{\zeta_i}{\To} \pi_1(S, q) \stackrel{\rho}{\To}
    PSL_2\R.
\]

Since each $\rho(C_i)$ is hyperbolic, and $\rho$ on each boundary curve of $S$ is non-elliptic, there is a well-defined relative Euler class $\E(\rho_i)$, and they are
additive by lemma \ref{euler_additivity}:
\[
    \sum_{i=1}^{-\chi(S)} \E(\rho_i)[S_i] = \E(\rho)[S] = - \chi(S).
\]
Since by the Milnor-Wood inequality \cite{Me10MScPaper0, Milnor, Wood} $|E(\rho_i)[S_i]| \leq 1$, we have $\E(\rho_i)[S_i] = 1$ for each $i$.

By proposition \ref{punctured_torus_euler_class} then, for each punctured torus $S_i$, $\Tr[g,h] \leq -2$ for some basis $g,h$ of $\pi_1(S_i, q_i)$; in particular, $\rho_i$ is not virtually abelian. Thus, by sections 5.2--5.3 of the prequel \cite{Me10MScPaper1} in the punctured torus case, and by section \ref{sec:hyperbolize_pants} above in the pants case, each $\rho_i$ is the holonomy of a complete hyperbolic structure on $S_i$ with totally geodesic boundary. 

It remains to fit the pieces together.

\subsection{Putting the pieces together}
\label{sec:putting_pieces_together}

We construct the hyperbolic structure
on $S$ piece by piece, starting from a first piece $S_1$ whose basepoint coincides with that of $S$, $q=q_1$. We then work outwards
along the tree $\T$ dual to the decomposition. By our choice of decomposition (figure \ref{fig:91}), when adding a new
piece, we only need to ensure it attaches along one boundary curve.

Consider first the case where $S$ is an $n$-holed sphere. This decomposes into $n-2$ pants. The combinatorial arrangement must be as in figure \ref{fig:91},
minus the punctured tori. We have a presentation $\pi_1(S) = \langle C_1, \ldots, C_n \; | \; C_1 \cdots C_n = 1 \rangle$.

Let $\rho(C_i) = c_i$, all non-elliptic, and let $\tilde{c}_i \in \widetilde{PSL_2\R}$ be preferred lifts. As $\E(\rho)[S] = -\chi(S) = n-2$, proposition \ref{euler_algebra} gives $\tilde{c}_1 \cdots \tilde{c}_{n} = \z^{n-2}$. Consider the following algebraic decomposition of the relator, corresponding to
the decomposition of the surface.
\begin{align*}
    1 &= \left[ C_1 C_2 (C_1 C_2)^{-1} \right] \left[ (C_1 C_2) C_3 (C_1 C_2 C_3)^{-1} \right] \cdots \\
    & \quad \cdots \left[ (C_1 C_2 \cdots
    C_{n-3}) C_{n-2} (C_1 \cdots C_{n-2})^{-1} \right]
    \left[ (C_1 C_2 \cdots C_{n-2}) C_{n-1} C_{n} \right].
\end{align*}
Each expression in square brackets is the relator in the presentation of each $\pi_1(S_i)$; as each $\E(\rho_i)[S_i] = 1$, this relator equals $\z$. From proposition \ref{pants_construction}, each $\rho_i$ is the holonomy of a complete structure on $S_i$ with each boundary curve $(C_1 \ldots C_i, C_{i+1}, (C_1 \cdots C_{i+1})^{-1})$ bounding $S_i$ on its right. Each decomposition curve is some $(C_1 C_2 \ldots C_j)$, and appears in two relators, which cancel. Hence in the corresponding fundamental domains, the curve corresponding to $(C_1 \ldots C_j)$ bounds one fundamental domain on its right, and the other on its left. Note that although our fundamental domains are degenerate along edges corresponding to parabolic boundary components, the decomposition curves are all hyperbolic, hence not degenerate.

To construct a hyperbolic structure, we piece together developing maps. Take the first surface $S_1$, with $q=q_1$, $\rho_{S_1} = \rho_1$. Take a preferred
lift $\tilde{q} = \tilde{q}_1$ of $q_1$, and partial lift of $\T$, in the universal cover $\tilde{S}_1$. Hyperbolizing our pants, we construct an octagonal fundamental domain, with corresponding basepoint $\bar{q} = \bar{q}_1 \in \hyp^2$, and extend to a developing map $\D_{S_1}: \tilde{S}_1 \To \hyp^2$. We then have the following, for $W=S_1$:
\begin{enumerate}
    \item
        A developing map $\D_W: \tilde{W} \To \hyp^2$ giving a hyperbolic structure on $W$ with totally geodesic or cusped boundary components and holonomy $\rho_W: \pi_1(W,q) \To PSL_2\R$, where $\rho_W$ is the restriction of $\rho$ to $\pi_1(W,q)$.
    \item
        Suppose $C_j$ is a boundary curve of $W$ which intersects $\T$, i.e. a decomposition curve. Via $\T$, take a representative of $C_j \in \pi_1(W, q)$ and take a canonical boundary edge $\tilde{C}_j$ of the universal cover $\tilde{W}$, where $\tilde{C}_j \cong \R$ covers $C_j \cong S^1$. Then $\D_W( \tilde{C}_j ) = \Axis c_j$, where $c_j = \rho_W(C_j) = \rho(C_j)$, .
\end{enumerate}

We inductively construct $\D_W$ verifying the above, for successively larger $W$. Suppose we have such a $W$; we adjoin a new pair of pants $S_k$ adjacent to $W$ and obtain a geometric structure on $W' = S_k \cup W$ with the same properties.

To do this, we use $\T$ to obtain preferred inclusions $\tilde{S}_k \hookrightarrow \tilde{W}' \hookrightarrow \tilde{S}$ and $\tilde{W} \hookrightarrow \tilde{W}' \hookrightarrow \tilde{S}$. Now $S_k \cap W$ is a single decomposition curve, say $C_k$, and using $\tilde{\T}$, $C_k$ has preferred lifts $\tilde{C}_k$ in $\tilde{S}_k$ and $\tilde{W}$, which agree upon inclusion into $\tilde{W}'$. So within $\tilde{W}'$, the preferred universal covers $\tilde{S}_k$ and $\tilde{W}$ intersect precisely along the preferred lift $\tilde{C}_k$.

The representation $\rho_k$, as described in section \ref{sec:hyperbolize_pants}, is the holonomy of a complete hyperbolic structure on $S_k$, with basepoint $q_k$ lifting to $\tilde{q}_k \in \tilde{S}_k \hookrightarrow \tilde{W}'$. We obtain a developing map $\D_k: \tilde{S}_k \To \hyp^2$ which takes $\tilde{C}_k$ to $\Axis c_k$. After possibly adjusting $\D_k$ by a diffeomorphism along $\tilde{C}_k$, $\D_k$ agrees with $\D_W$ along $\tilde{C}_k$. Combining the two developing maps $\D_k, \D_W$ gives a partial developing map of $W'$, which extends equivariantly to a true developing map $\D_{W'}$ for $W'$, satisfying the above conditions.

\begin{figure}
\centering
\includegraphics[scale=0.3]{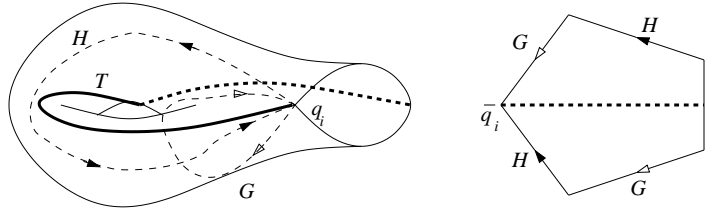}
\caption{Choice of $\T$ on each punctured torus, and developing
image.} \label{fig:88a}
\end{figure}

Continuing in this manner, we obtain a hyperbolic structure on the entire surface $S$, where $S$ has genus $0$. Attaching punctured tori is no more difficult: the same argument applies. In order for our constructions and basepoints to match, we choose the dual tree $\T$ to avoid the curves $G_k, H_k$ of $\pi_1(S_k, q_k)$ forming the basis for the construction: see figure \ref{fig:88a}. 

We thus construct a developing map for $S$, giving a geometric structure with holonomy $\rho$. The boundary components are as desired: a parabolic boundary must lie in one of the pants of our decomposition, and the construction of section \ref{sec:hyperbolize_pants} gives a cusp. This concludes the proof of Goldman's theorem.

\section{Constructions for the genus 2 surface}
\label{sec:constructions_for_genus_2}

We prove theorem \ref{restricted_version}. Let $S$ be a closed genus 2 surface; assume $\E(\rho)[S] = 1$; if $\E(\rho)[S] = -1$ then the same arguments apply with opposite orientation. We suppose that there is a separating curve $C$ on $S$ such that $\rho(C)$ is not hyperbolic.

\subsection{Splitting into tori}

\begin{figure}
\centering
\includegraphics[scale=0.3]{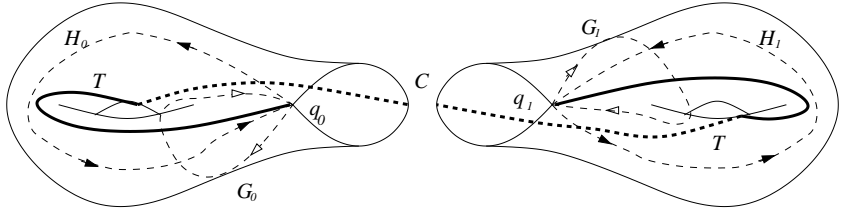}
\caption{Decomposition of $S$ and dual tree $\T$.} \label{fig:85a}
\end{figure}

\begin{figure}
\centering
\includegraphics[scale=0.3]{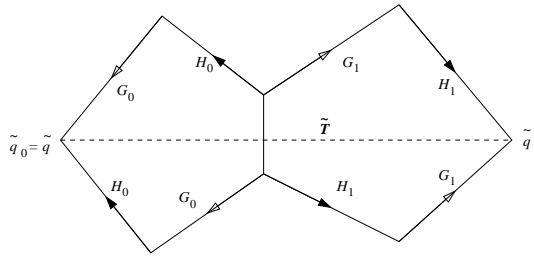}
\caption{Preferred lifts of basepoints and $\T$ in the universal
cover $\tilde{S}$.} \label{fig:85b}
\end{figure}

Let $q$ be a basepoint on $C$; let $C$ split $S$ into two punctured tori $S_0, S_1$. A dual tree to the splitting is just an edge with a vertex at either end. We take basepoints $q_0=q_1=q$ for $S_0, S_1, S$ respectively. On $S_i$, let $G_i, H_i \in \pi_1(S_i, p_i)$ be basis curves, so $[G_0, H_0]$ and $[G_1, H_1]$ are homotopic to $C$, but traversed in opposite directions. Choose the dual tree $\T$ to run between $G_i, H_i$ as in figure \ref{fig:88a}. Take preferred lifts $\tilde{T}$, $\tilde{q} = \tilde{q}_0$ and $\tilde{q}_1 \in \tilde{S}$; as in section \ref{sec:splitting} we have homomorphisms
\begin{align*}
    &\iota_0: \pi_1 (S_0, q_0) \hookrightarrow \pi_1(S, q_0) = \pi_1(S,q), \quad
    &\iota_1: \pi_1 (S_1, q_1) \hookrightarrow \pi_1(S, q_1)= \pi_1(S,q),
    \\
    &\zeta_0: \pi_1(S, q_0) \stackrel{\text{Id}}{\To} \pi_1(S, q),
    \quad
    &\zeta_1: \pi_1(S, q_1) \stackrel{\cong}{\To} \pi_1(S, q),
\end{align*}
and representations $\rho_0 = \rho \circ \zeta_0 \circ \iota_0$,
$\rho_1 = \rho \circ \zeta_1 \circ \iota_1$. Note by our choice of
$\T$, we have $\tilde{q}_0 \neq \tilde{q}_1$, even though $q_0 =
q_1$. See figures \ref{fig:85a} and \ref{fig:85b}.

With these choices, $\pi_1(S,q) = \langle G_0, H_0, G_1, H_1 \; | \; [G_0, H_0]\; [G_1, H_1] = 1 \rangle$. Let $L \in \pi_1(S,q)$ denote the loop traced out by $\T$ from $q_0$ to $q_1$; $\zeta_1$ is conjugation by $L$; in fact $L = G_0^{-1} H_0^{-1} G_1 H_1 = H_0^{-1} G_0^{-1} H_1 G_1$. With $g_i = \rho(G_i)$, $h_i = \rho(H_i)$, we have
\[
    \rho_0 (G_0) = g_0, \quad \rho_0 (H_0) = h_0, \quad \rho_1(G_1)
    = \rho(L) g_1 \rho(L^{-1}), \quad \rho_1(H_1) = \rho(L) h_1
    \rho(L^{-1}).
\]

Since $\E(\rho)[S] = 1$ we have by proposition \ref{euler_algebra} $[g_0, h_0][g_1, h_1] = \z$. Note that $[g_i, h_i]$, $[g_i^{-1}, h_i^{-1}]$ and $[\rho_i(G_i)^{-1}, \rho_i(H_i)^{-1}]$ are all conjugates in the holonomy group, hence lie in the same region of $\widetilde{PSL_2\R}$. In fact, choosing arbitrary lifts $\tilde{g}_0, \tilde{h}_0, \tilde{g}_1, \tilde{h}_1 \in \widetilde{PSL_2\R}$ and $\widetilde{\rho(L)} = \tilde{g}_0^{-1} \tilde{h}_0^{-1} \tilde{g}_1 \tilde{h}_1 = \z \tilde{h}_0^{-1} \tilde{g}_0^{-1} \tilde{h}_1 \tilde{g}_1$, we can easily obtain $\rho_0 ([G_0^{-1}, H_0^{-1}]) \rho_1 ([G_1^{-1}, H_1^{-1}]) = \z$.

As $[G_i, H_i]$ is homotopic to $C$, traversed in some direction, each $[g_i, h_i]$ is not hyperbolic. Since we know the possible regions of $\widetilde{PSL_2\R}$ in which commutators lie (e.g. \cite{Me10MScPaper0, Me10MScPaper1, Milnor, Wood, Eisenbud-Hirsch-Neumann, Goldman88}) we have
\[
    [g_0, h_0], [g_1, h_1] \in \{1\} \cup
    \Ell_{-1} \cup \Ell_1 \cup \Par_{-1}^+ \cup \Par_0 \cup
    \Par_1^-.
\]
As $[g_0, h_0], [g_1, h_1]$ are inverses in $PSL_2\R$, they are both elliptic, both parabolic, or both the identity. Applying properties of $\Theta$ (see \cite{Me10MScPaper0}) we have $\Theta([g_0, h_0]) + \Theta([g_1, h_1]) = \pi$. We know $\Theta$ of the various regions of $\widetilde{PSL_2\R}$; assuming without loss of generality $\Theta([g_0, h_0]) \leq \Theta([g_1, h_1])$, there are only the following two possibilities. (In particular, neither of $[g_0, h_0], [g_1, h_1]$, considered in $PSL_2\R$, can be the identity.) 
\begin{enumerate}
    \item
        {\bf Elliptic case.} $[g_0, h_0] \in \Ell_1$ with
        $\Theta([g_0, h_0]) \in \left( 0, \pi/2 \right]$, and $[g_1,
        h_1] \in \Ell_1$ with $\Theta([g_1, h_1]) \in \left[ \pi/2,
        \pi \right)$.
    \item
        {\bf Parabolic case.} $[g_0, h_0] \in \Par_0^+$ and $[g_1,
        h_1] \in \Par_1^-$.
\end{enumerate}
We will consider these two cases separately in the next two sections.

\subsection{Piecing together along an elliptic}

We have $[g_0, h_0], [g_1, h_1] \in \Ell_1$ with $\Theta([g_0, h_0]) \in (0, \pi/2]$ and $\Theta([g_1, h_1]) \in [\pi/2, \pi)$. Applying proposition 5.3 of \cite{Me10MScPaper1}, $\rho_0$ is the holonomy of a hyperbolic cone-manifold structure on $S_0$ with corner angle in $[2\pi, 3\pi)$, a ``large angle" elliptic case; and $\rho_1$ is a ``small angle" elliptic case, with corner angle in $(\pi, 2\pi]$.

Recall that this construction (section 5.4 of \cite{Me10MScPaper1}) gives a pentagonal fundamental domain $\Pent(g,h;p)$, choosing $p$ close to $\Fix [g^{-1},h^{-1}]$ judiciously. Since $[g_0^{-1}, h_0^{-1}]$ and $[g_1^{-1}, h_1^{-1}]$ are both in $\Ell_1$, $[G_0, H_0]$ will bound $S_0$ on its left and $[G_1, H_1]$ will bound $S_1$ on its left; with preferred lifts $\tilde{q}_0$, $\tilde{q}_1$ and $\tilde{T}$ as above, the points $\tilde{q}_0$ and $\tilde{q}_1$ will lie on adjacent pentagonal fundamental domains as in figure \ref{fig:85b}.

We start with a fundamental domain $\Pent(g_0, h_0; p_0)$ for $S_0$, and add on the fundamental domain for $S_1$. The representation $\rho_1: \pi_1(S_1, q_1) \to PSL_2\R$ is given by $\rho \circ \zeta_1 \circ \iota_1$ and corresponds to the holonomy of a developing map where $q_1$ lifts to $\tilde{q}_1 \in \tilde{S}_1 \hookrightarrow \tilde{S}$. We will construct $\Pent(\rho_1(G_1), \rho_1(H_1); p_1)$ where $p_1$ is close to $\Fix [\rho_1(G_1)^{-1}, \rho_1(H_1)^{-1}]$, which is the same as the fixed point of its inverse $[\rho_0(G_0)^{-1}, \rho_0(H_0)^{-1}] = [g_0^{-1}, h_0^{-1}]$.

\begin{figure}
\centering
\includegraphics[scale=0.3]{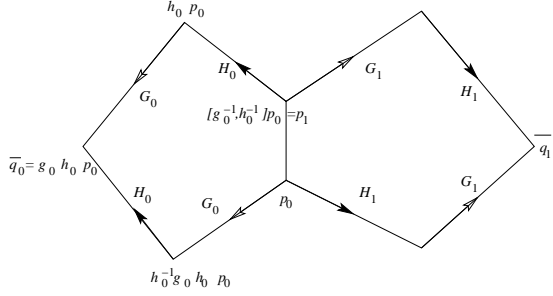}
\caption{Putting the pieces together.} \label{fig:84}
\end{figure}

We must find basepoints $p_0, p_1$ such that the pentagons $\Pent(\rho_0(G_0), \rho_0(H_0); p_0)$ and $\Pent(\rho_1(G_1), \rho_1(H_1); p_1)$ join precisely along the edges representing their boundary, without folding. We must have $p_1 = [g_0^{-1}, h_0^{-1}] p_0$; see figure \ref{fig:84}.

Recall now proposition 5.3 of \cite{Me10MScPaper1} in detail. Let $r = \Fix [g_0^{-1}, h_0^{-1}]$. There exists a closed semicircular disc $C_{\epsilon_0}(r)$ centred at $r$ such that if $p_0 \in C_{\epsilon_0}(r)$, $p_0 \neq r$, then $\Pent(\rho_0(G_0),\rho_0(H_0);p_0)$ is non-degenerate and bounds an
embedded disc. Similarly we have $C_{\epsilon_1}(r)$ about $p_1$ for which $\Pent(\rho_1(G_1), \rho_1(H_1); p_1)$ bounds an embedded disc. Take $\epsilon =
\min(\epsilon_1, \epsilon_2)$. On the circle of radius $\epsilon$ about $r$, there is a closed arc of angle $\pi$ on which $p_0$ can validly lie; and a closed arc of angle $\pi$ on which $[g_1^{-1},h_1^{-1}]p_1$ can validly lie. Hence these arcs must overlap, and we may take $p_0$ and $p_1$ so that 
$\Pent(\rho_0(G_0),\rho_0(H_0);p_0)$ and $\Pent(\rho_1(G_1),\rho_1(H_1);p_1)$ are both non-degenerate, bounding embedded discs, and matching along their boundary arcs.

We have an immersed non-degenerate geodesic octagon bounding an immersed disc in $\hyp^2$, which is a fundamental domain for $S$; we may then extend equivariantly to a true developing map, obtaining a cone-manifold structure on $S$, with holonomy $\rho$. There is at most one cone point, at the vertices (all of which are identified) of the fundamental domain. The cone angle is the sum of the interior angles of the octagon, which is equal the sum of the two corner angles $\theta_i$ in the punctured tori $S_i$. From proposition 5.3 (or lemma 3.6) of \cite{Me10MScPaper1} again, $\theta_i = 3\pi - \Twist([g^{-1}, h^{-1}], p_i)$. Since $\rho_0 ([G_0^{-1}, H_0^{-1}]) \rho_1 ([G_1^{-1}, H_1^{-1}]) = \z$, $\Twist([\rho_0(G_0),\rho_0(H_0)],p_0) + \Twist([\rho_1(G_1), \rho_1(H_1)],p_1) = 2\pi$. Hence the cone angle is $6\pi - 2\pi = 4\pi$.

So we obtain the desired hyperbolic cone-manifold structure on $S$.

\subsection{Piecing together along a parabolic}
\label{sec:piecing_along_parabolic}

We have $[g_0, h_0] \in \Par_0^+$ and $[g_1, h_1] \in \Par_1^-$; since we know the traces of various regions (e.g. lemma 2.8 of \cite{Me10MScPaper1}), $\Tr[g_0, h_0] =2$, $\Tr[g_1, h_1] =-2$. Hence we may apply the corresponding results (sections 5.5 and 5.3) of the prequel \cite{Me10MScPaper1}; the strategy is the similar to the previous section. Let $r = \Fix [\rho_0(G_0)^{-1}, \rho_0(H_0)^{-1}] = \Fix [\rho_1(G_1)^{-1}, \rho_1(H_1)^{-1}]$.

First consider $S_0$. From section 5.5 of \cite{Me10MScPaper1} we may take a basis $G_0, H_0$ of $\pi_1(S_0)$ and a point $p_0$ arbitrarily close to $r$ such that $\Pent(\rho_0(G_0), \rho_0(H_0); p_0)$ is non-degenerate and bounds an embedded disc. Since $[\rho_0(G_0), \rho_0(H_0)] = [g_0, h_0] \in \Par_0^+$, $\partial S_0$ traversed in the direction of $[G_0, H_0]$ bounds $S_0$ on its left. This gives a hyperbolic cone-manifold structure on $S_0$ corresponding to a preferred lift $\tilde{q}_0$ of $q_0$ with holonomy $\rho_0$. The corner angle is $3\pi - \Twist([\rho_0(G_0)^{-1}, \rho_0(H_0)^{-1}],p_0)$.

Now consider $S_1$. As above, we take a basis $G_1, H_1$ of $\pi_1(S_1, q_1) \stackrel{\zeta_1 \circ \iota_1}{\To} \pi_1(S,q)$ such that $[G_0, H_0][G_1, H_1] = 1$. From section 5.3 of \cite{Me10MScPaper1}, the representation $\rho_1$ is discrete, and the quotient of $\hyp^2$ by the image of $\rho_1$ is a cusped torus. We may take
$p_1$ anywhere sufficiently close to $r$, and obtain $\Pent(\rho_1(G_1),\rho_1(H_1);p_1)$ non-degenerate bounding an embedded disc. This gives a hyperbolic cone-manifold structure on $S_1$ corresponding to the preferred lift $\tilde{q}_1$ of $q_1$ with holonomy $\rho_1$. The corner angle is $3\pi -\Twist([\rho_1(G_1)^{-1}, \rho_1(H_1)^{-1}],p_1)$, and boundary $\partial S_1$ traversed in the direction of $[G_1, H_1]$ bounding $S_ 1$ on its left.

Hence we may take $p_0,p_1$ such that $p_1 = [g_0^{-1},h_0^{-1}]p_0$ and both $\Pent(\rho_0(G_0),\rho_0(H_0);p_0)$, $\Pent(\rho_1(G_1),\rho_1(H_1);p_1)$ are non-degenerate pentagons bounding immersed discs. Since they both have the same orientation, they fit together without folding along their boundary edges to give
a non-degenerate octagon in $\hyp^2$ bounding an immersed disc, and hence a cone-manifold structure on $S$ with holonomy $\rho$. Since $[\rho_0(G_0)^{-1}, \rho_0(H_0)^{-1}] [\rho_1(G_1)^{-1},\rho_1(H_1)^{-1}] = \z$ we have $\Twist([\rho_0(G_0)^{-1},\rho_0(H_0)^{-1}],p_0) + \Twist([\rho_1(G_1)^{-1}, \rho_1(H_1)^{-1}],p_1) = 2\pi$. Hence the cone angle is $6\pi - 2\pi = 4\pi$.

Geometrically, one half of $S$ has the nice structure of a truncated
cusped torus, and the other half is a rather uglier handle tacked on
to the truncated cusp. This concludes the proof of theorem \ref{restricted_version}.

\section{Representations with $\E(\rho)[S] = \pm \left( \chi(S) + 1 \right)$}
\label{sec:one_off_extremal}

In this section we prove theorem \ref{technical_prop}, or its more precise version theorem \ref{thm:technical_precise}. So let $S$ be a closed surface of genus $g \geq 2$ and $\rho$ be a representation with $\E(\rho)[S] = \pm \left( \chi(S) + 1 \right)$; we assume there exists a non-separating simple closed curve $C$
with $\rho(C)$ elliptic. We only consider the case $\E(\rho)[S] = - \chi(S) - 1$; the case $\chi(S) + 1$ is simply orientation reversed.

Note that theorem \ref{technical_prop} is not vacuous: such representations do exist. For instance, theorem 4.1 of \cite{Me10MScPaper1} implies that there exist representations for a punctured torus which take a non-separating simple closed curve to an elliptic and have $\Tr[g,h]>2$; glue this with the holonomy representation of a complete hyperbolic structure on a surface of genus $g-1$ with one boundary component.

\subsection{Easier case: gluing along a parabolic}
\label{sec:easier_parabolic}
\label{sec:simple_cases}

Given non-separating simple $C$, we can find a separating simple closed curve $D$, disjoint from $C$, cutting $S$ into two
pieces: a punctured torus $S_1$ containing $C$; and a surface $W$ of genus $g-1$ with $1$ boundary component. Choosing basepoints $q = q_W, q_1$ for $S, W$ and $S_1$ respectively, and a dual tree $\T$ as in section \ref{sec:splitting}, we obtain representations $\rho_W, \rho_1$ on $W$ and $S_1$.

Take a basis $G,H$ for $\pi_1 (S_1, q_1)$ with $G$ freely homotopic to $C$, and $[G,H]$ homotopic to $D$, so $\rho_1(G)=g$ is elliptic. Lemma 2.2 of \cite{Me10MScPaper1} (following \cite{Goldman03}, see also \cite{Me10MScPaper0}) implies that if $\Tr[g,h] < 2$ then $g,h$ are hyperbolic; so $\Tr[g,h] \geq 2$. Hence $[g,h]$ is not elliptic, and relative Euler classes are well-defined. By proposition \ref{punctured_torus_euler_class}, $\E(\rho_1)[S_1] = 0$; hence by additivity of the relative Euler class \ref{euler_additivity}, $\E(\rho_W)[W] = - \chi(W)$. Thus Goldman's theorem applies to $\rho_W$, which is the holonomy of a complete hyperbolic
structure on $W$ with totally geodesic or cusped boundary. The holonomy of the boundary of $W$ must therefore be parabolic or hyperbolic.

The parabolic case is easier, and we deal with it first. Having $g$ elliptic, $[g,h]$ parabolic, and $\Tr[g,h]=2$, implies that $\rho_1$ is reducible, not abelian, indeed not virtually abelian. So by proposition 5.5 of \cite{Me10MScPaper1}, $\rho_1$ is the holonomy of a hyperbolic cone-manifold structure on $S_1$ with one corner point, with angle $>2\pi$. The same method as in section \ref{sec:piecing_along_parabolic} then allows us to glue together two developing maps for $W$ and $S_1$, and gives the following.

\begin{lem}
\label{simple_case}
    If $\rho(D)$ is parabolic then $\rho$ is the holonomy of a hyperbolic cone-manifold structure on $S$ with one cone point of
    angle $4\pi$.
\qed
\end{lem}

This leaves only with the case where $\rho(D)$ is hyperbolic, i.e. $\Tr [g,h] > 2$.

\subsection{Piecing together along a hyperbolic}

In this case, $\rho_1$ has $g$ elliptic, $\Tr[g,h]=t>2$ and $\E(\rho_1)[S_1]=0$; $\rho_W$ is the holonomy of a complete hyperbolic structure on $W$, hence discrete.

The quotient of $\hyp^2$ by the image of $\rho_W$ is a flared surface. As discussed in section 6 of \cite{Me10MScPaper1} and section \ref{sec:piecing_along_parabolic} above, we can truncate the ``flares" along geodesics in the homotopy classes of the boundary curves; we can also truncate a ``flare" away from the geodesic, and obtain a piecewise geodesic boundary, with a single corner point. This truncation can be done arbitrarily outside the convex core, obtaining corner angles in $(0, \pi)$. It may also be done inside the convex core, producing a corner angle in $(\pi, 2\pi)$; but we cannot truncate too far inside the surface. Nonetheless if we stay within the collar width of the geodesic then we are guaranteed still to obtain a cone-manifold structure on $W$. The collar width $w(t)$ only depends on $t$ and is given by $\sinh w(t) = 1/ \sinh \left( \frac{d(t)}{2} \right)$, where $d(t) = 2 \cosh^{-1}(t/2)$ is the length of the geodesic boundary curve homotopic to $D$. See \cite{Buser} for details.

We wish to perform such a truncation inside the convex core of $W$, to find a hyperbolic cone-manifold structure on $W$ with a corner angle in $(\pi, 2\pi)$, which pieces together with a cone-manifold structure on $S_1$ with corner angle in $(2\pi, 3\pi)$, to give a cone-manifold structure on $S$ with a single cone point of angle $4\pi$.  In the parabolic cases of sections \ref{sec:piecing_along_parabolic} and \ref{sec:easier_parabolic} above, the convex core extends to infinity, and the non-complete half of the representation can be constructed arbitrarily close to infinity, so developing maps can easily be pieced together. But the constructions of section 5.6 of \cite{Me10MScPaper1} do not work arbitrarily close to infinity; the hyperbolic case is more difficult.

We will find a pentagonal fundamental domain for $S_1$, which pieces together with the developing map $\D_W$ along (a lift of) $D$, to obtain a developing map for a hyperbolic cone-manifold structure on $S$. As in the parabolic case, there is no folding. The corner angle on $S_1$ is $3\pi - \Twist(\rho(D),p)$, by section 5.6 of \cite{Me10MScPaper1}. The corner angle on $W$ is $\pi + \Twist(\rho(D),p)$. Here we take the simplest lift of $\rho(D)$ into $\widetilde{PSL_2\R}$. So as in the parabolic case, once the developing map is constructed the cone angle is automatically $4\pi$.

These considerations essentially reduce the problem to two measure-theoretic propositions, to which the next two sections are dedicated.

To understand these propositions, we make some preliminary remarks. Consider a punctured torus $S_1$, a basis $G,H$ for $\pi_1(S_1)$, and a representation $\rho: \pi_1(S_1) \To PSL_2\R$ with $\Tr[g,h]>2$. Define $\rho$ to be \emph{$\epsilon$-good} for a specified orientation of $S$, if there exists a basis $(G',H')$, of the same orientation as $(G,H)$, and a point $p$ at distance $\leq \epsilon$ from $\Axis [g'^{-1},h'^{-1}]$ (where $g' = \rho(G')$, $h' = \rho(H')$), such that the pentagon $\Pent(g',h';p)$ is non-degenerate, bounds an embedded disc, and is of the specified orientation. That is, $\epsilon$-good representations give cone manifold structures on punctured tori with one corner point, of specified orientation, with a pentagonal fundamental domain having boundary edge within $\epsilon$ of the axis of the boundary holonomy. (Who would say this is a bad thing?) Note that if a representation is $\epsilon$-good, so is any conjugate representation.

A character is \emph{$\epsilon$-good} for a specified orientation of $S$ if it is the character of an $\epsilon$-good representation for the same orientation. Since we are only concerned with $\Tr[g,h]>2$, by proposition 4.2 of \cite{Me10MScPaper1}, all representations concerned are irreducible; and hence characters correspond precisely to conjugacy classes of representations. So a character is $\epsilon$-good iff one corresponding representation is $\epsilon$-good, iff all corresponding representations are $\epsilon$-good.

Define characters or representations which are not $\epsilon$-good to be \emph{$\epsilon$-bad}. 

Recall from section \ref{sec:character_variety} above the relative character variety $X_t(S_1) = X(S_1) \cap \kappa^{-1}(t)$. Recall from section \ref{sec:action_on_character_variety} that the measure $\mu_t$ on $X_t(S_1)$ is invariant under the action of $\Gamma$. We are considering $\rho$ which take some simple closed curve to an elliptic. So let $\Omega_t \subset X_t(S_1)$ be the set of characters of representations taking some simple closed curve to an elliptic.

For a specified orientation of $S_1$, let $B_t$ be the set of $w(t)$-bad characters in $\Omega_t \subset X_t(S_1)$, where $w(t)$ is the collar width defined above.

\begin{prop}
\label{ergodicity_prop}
    For all $t>2$, $\mu_t(B_t) = 0$. That is, $\mu_t$-almost every character in $\Omega_t$ is $w(t)$-good, for the specified orientation of $S_1$. 
\end{prop}

The proof of this result will use ergodicity properties of action of $\Gamma$ on the character variety. The strategy is to show that \emph{some} representation produces a desirable pentagon; and to use ergodicity to show that changing basis we can ``almost" move anywhere within the character variety, so we can get close to $\rho$ and produce such a pentagon.

Note the word ``almost" cannot be removed from the statement of \ref{ergodicity_prop}: characters of virtually abelian representations lie inside $\Omega_t \subset X_t(S_1)$, and such representations are $\epsilon$-bad for any $\epsilon$.

Let $\mathcal{U}$ denote the set of separating curves $D$ which split $S$ into a punctured torus $S_1$ and another surface $W$. For $D \in \mathcal{U}$ and $t > 2$, let $B_{D,t} \subset X^\pm (S)$ denote the set of all characters of (twisted) representations with euler class $-\chi(S)-1$, which take $D$ to have trace $t$, and which restrict on $S_1$ to a character in $B_t$. Such a restriction is well-defined, since the trace depends only on the conjugacy class, hence free homotopy class, of each loop.

In particular, a character in $B_{D,t}$ is non-abelian; the restriction to $S_1$ corresponds precisely to a conjugacy class of representations; such a representation $\rho_1$ takes some simple closed curve on $S_1$ to an elliptic; hence $\E(\rho_1)[S_1] = 0$ and $\E(\rho_W)[W]=-\chi(W)$; and with respect to any dual tree $\T$, the induced representation $\rho_1$ is $w(t)$-bad.

Let $B_D = \cup_{t>2} B_{D,t}$ and $B = \cup_D B_D$. So $B \subset \Omega_S \subset X^\pm(S)$; recall $\Omega_S$, as in the statement of theorem \ref{thm:technical_precise}, denotes characters of representations $\rho$ with $\E(\rho)[S] = \pm (\chi(S) + 1)$, which take some simple closed curve to an elliptic.

\begin{prop}
\label{measure_theory_prop}
    $\mu_S (B) = 0$.
\end{prop}

\begin{Proof}[of theorem \ref{technical_prop}, precisely stated as \ref{thm:technical_precise}, assuming \ref{measure_theory_prop}]
    By proposition \ref{measure_theory_prop}, it suffices to show that a representation $\rho$ with character in $\Omega_S \backslash B$ is the holonomy of a cone-manifold structure on $S$ with a single cone point with cone angle $4\pi$. Such a $\rho$ sends sends some simple closed $C$ to an elliptic, and without loss of generality $\E(\rho)[S] = - \chi(S) - 1$. Recall section \ref{sec:easier_parabolic}. For any separating $D$ cutting off a punctured torus $S_1$ containing $C$, $\E(\rho_1)[S_1] = 0$, $\rho_W$ is discrete, and $\rho(D)$ is parabolic or hyperbolic. If it is parabolic, by lemma \ref{simple_case} we are done. So we may assume that $\rho(D)$ is hyperbolic; choosing a basepoint $q_1 \in \partial S_1$ for $S_1$ and basis $G,H$ for $\pi_1(S_1)$, we have $\Tr[g,h]=t>2$, i.e. $\rho_1$ has character in $X_t(S_1)$ for $t>2$.  The presence of elliptic $\rho(C)$ means in fact the character of $\rho_1$ lies in $\Omega_t$.

    As $\rho$ has character in $\Omega_S \backslash B$, $\rho_1$ has character in $\Omega_t \backslash B_t$, i.e. $\rho_1$ is $w(t)$-good. So we may take a basis $G',H'$ of $\pi_1(S,q_1)$, of the same orientation, and $p$ within the collar width of $\Axis [g'^{-1}, h'^{-1}]$, such that $\Pent(g',h';p)$ is non-degenerate, bounds an embedded disc, and the edge $p \to [g'^{-1}, h'^{-1}]p$ bounds the pentagon on its left.

    As $\rho_W$ is discrete, $W$ is the quotient of a convex core. By truncating $W$ within its collar, and possibly applying a homeomorphism supported near $\partial W$, as in previous cases we may piece together a partial developing map and extend to a developing map on $\tilde{S}$ giving a hyperbolic cone-manifold structure on $S$ with a single cone point. There is no folding, since we chose the orientation on $S_1$ to avoid it. As in previous cases, the cone angle is automatically $4\pi$. So we have the desired hyperbolic cone-manifold structure.
\end{Proof}

\subsection{Ergodicity}

Let $\rho: \pi_1(S_1) \To PSL_2\R$ be a representation and let $G,H$ be a basis of $\pi_1(S_1)$. Lift $g,h$ to $SL_2\R$ arbitrarily and let $(\Tr g, \Tr h, \Tr gh) = (x,y,z) \in X(S_1)$. Recall section \ref{sec:action_on_character_variety} and consider the action of $\mcg (S) \cong \Out \pi_1(S_1)$ and $\Gamma \cong PGL_2\Z \ltimes \left( \frac{\Z}{2} \oplus \frac{\Z}{2} \right)$ on $X(S_1)$; this action preserves the level sets $X_t(S_1) = \kappa^{-1}(t) \cap X(S_1)$ and the symplectic form and measure $\mu_t$ on each $X_t(S_1)$. For the proof of \ref{ergodicity_prop} we are only interested in $t>2$, for which $X_t(S_1) = \kappa^{-1}(t)$ (theorem 4.1 of \cite{Me10MScPaper1}).

In the case $t > 2$ there are no reducible representations and a character correpsonds uniquely to a conjugacy class of representations. Goldman \cite{Goldman03} proved that $X_t(S_1)$ consists of two types of representations:

\begin{enumerate}
    \item
        \emph{Pants representations:} Those $(x,y,z) \in X_t(S_1)$ equivalent to triples $(x',y',z')$ where $x',y',z' \leq -2$. Section \ref{sec:hyperbolize_pants} shows that these are discrete representations which can be considered the holonomy of a complete hyperbolic structure on a pair of pants with totally geodesic or cusped
        boundary. (Note that a given basis will not usually correspond to the boundary components of the pants.) Thus there are no elliptic elements in the image of $\rho$; and for any $(x',y',z') \sim (x,y,z)$ we have $|x'|, |y'|, |z'| \geq 2$.
    \item
        \emph{Representations with elliptics:} Those $(x,y,z) \in X_t(S_1)$
        equivalent to $(x',y',z')$ with some coordinate in $(-2,2)$.
        That is, there is some simple closed curve on $S_1$ with
        elliptic image: we have denoted these $\Omega_t$.
\end{enumerate}
Goldman gives an algorithm to change basis and reduce traces until
they are small or all negative --- essentially a greedy algorithm. Note
the action of $\Gamma$, or of $\Out \pi_1(S)$, preserves $\Omega_t$.

\begin{thm}[Goldman \cite{Goldman03}]
\label{thm:Goldman_ergodicity}
    For $t>2$, the action of $\Gamma$ on $\Omega_t$ is ergodic.
    \qed
\end{thm}
Recall \emph{ergodic} means that the only invariant sets in $\Omega_t$ under the action of $\Gamma$ are null or conull, i.e. they have measure zero, or their complement has measure zero.

The following result guarantees us good representations.
\begin{lem}
\label{good_rep_exists}
    For any $t>2$, $\epsilon>0$ and specified orienation of $S_1$, there exists an $\epsilon$-good representation $\rho_t^\epsilon$ for the specified orientation, with character in $\Omega_t$.
\end{lem}

\begin{figure}
\centering
\includegraphics[scale=0.4]{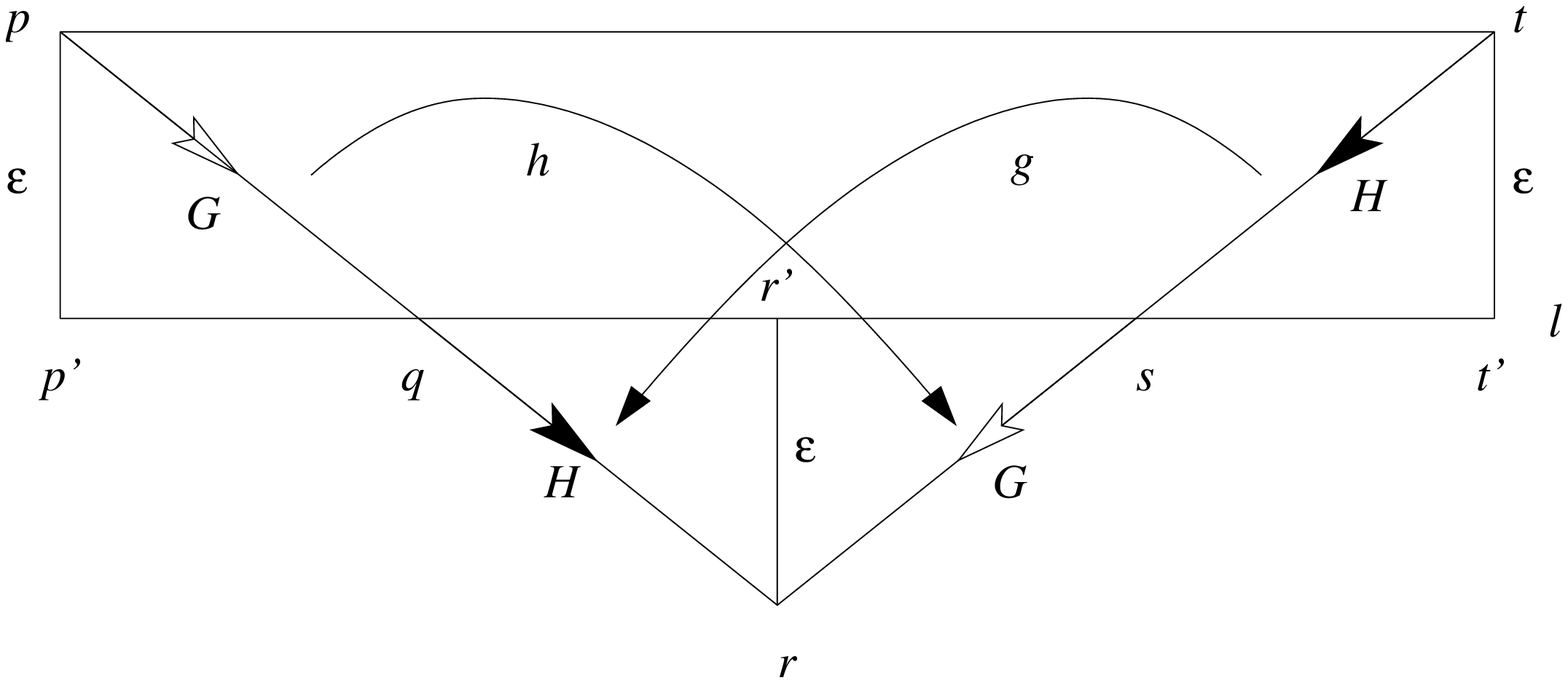}
\caption{Construction of a good representation.} \label{fig:88}
\end{figure}

\begin{Proof}
    Let $l \subset \hyp^2$ be a hyperbolic line and let $d = 2 \cosh^{-1} \left( t/2 \right)$, so $d$ will be the
    translation distance of $[g,h] \sim [g^{-1},h^{-1}]$. Define five points $p,q,r,s,t \in \hyp^2$, via Fermi coordinates on $l$: $p = (-d/2, \epsilon)$, $q = (-d/4, 0)$, $r = (0, -\epsilon)$, $s = (d/4, 0)$, $t = (d/2, \epsilon)$. Let $p' = (-d/2, 0)$, $r' = (0,0)$, $t'=(d/2,0)$ respectively
    be the projections of $p,r,t$ onto the line $l$, so the triangles $pp'q, rr'q, rr's, tt's$ are clearly congruent.
    It's clear that $pqr$ and $rst$ are geodesics. See figure \ref{fig:88}.

    Let $g$ be the orientation-preserving isometry carrying the directed segment $ts$ onto $qr$. Let $h$ be the orientation-preserving isometry carrying
    $pq$ onto $sr$. Then $q = h^{-1}ghp$, $r = ghp$, $s = hp$, $t = [g^{-1},h^{-1}] p$ so the pentagon $pqrst$ is actually $\Pent(g,h;p)$, non-degenerate and bounding an embedded disc. By replacing $\epsilon$ with $-\epsilon$ if necessary our pentagon has the desired orientation. This defines a representation $\langle G,H \rangle \To PSL_2\R$, which is clearly $\epsilon$-good. Note that replacing $\epsilon$ with any smaller $\epsilon'$ also gives an $\epsilon$-good representation.

    As $\epsilon \To 0$, the distance of $p$ from $l$ tends to $0$, $g$ and $h$ tend to half-turns, i.e. with trace $0$; so possibly replacing $\epsilon$ with a smaller $\epsilon'$, $g$ and $h$ are elliptic.

    Note that $gh$ translates along $pr$, from $p$ to $r$; and $g^{-1}h^{-1}$ translates along $tr$, from $r$ to $t$. Consider the action of $D(gh)$, then $D(g^{-1} h^{-1})$, on a unit tangent vector at $p$ pointing towards $p'$, along the perpendicular from $p$ to $l$. It ends up at $t$, pointing towards $t'$, along the perpendicular from $t$ to $l$. But this is precisely the action fo a translation along $l$ from $p$ to $t$, i.e. by distance $d$. It follows that $|\Tr[g,h]|=t>2$; but as $g,h$ are elliptic, $\Tr[g,h] < -2$ contradicts lemma 2.2 of \cite{Me10MScPaper1} (also \cite{Goldman03, Me10MScPaper0}); hence $\Tr[g,h]=t$. So we obtain an $\epsilon$-good $\rho_t^\epsilon$ for the specified orientation on $S_1$, with character in $\Omega_t$.
\end{Proof}

Note from the symmetry of our construction that $g$ and $h$ are conjugate via a reflection in a vertical axis of symmetry of figure \ref{fig:88}. In particular, (choosing lifts of elements in $PSL_2\R$ to $SL_2\R$ appropriately) $\Tr g = \Tr h$, and so the character is of the form $(x,y,z) = (x,x,z)$.

Now we can complete the proof of proposition \ref{ergodicity_prop}. 
\begin{Proof}[of \ref{ergodicity_prop}]
We must show that $\mu_t$-almost every character in $\Omega_t$ is $w(t)$-good, for a specified orientation of $S_1$. From lemma \ref{good_rep_exists} we can obtain a representation $\rho_t^{w(t)/2}$ in $\Omega_t$ which is $w(t)/2$-good (hence $w(t)$-good) for the specified orientation, with character of the form $(x^*, x^*, z^*)$ where $|x^*|<2$. Now any point $(x,y,z)$ in $\R^3$ sufficiently close (in the Euclidean metric) to $(x^*, x^*, z^*)$ on the surface $\kappa^{-1}(t)$ is also in $\Omega_t$ (having a coordinate with magnitude less than $2$), and is also $w(t)$-good. To see this, we can take $\rho$ with character $(x,y,z)$ and with $g,h$ close to those from $\rho_t^{w(t)/2}$; we can take $p$ close to that from $\rho_t^{w(t)/2}$, hence within $w(t)$ of $\Axis [g^{-1},h^{-1}]$; then the vertices of $\Pent(g,h;p)$ are close to those from $\rho_t^{w(t)/2}$. So the pentagon has the same orientation and still bounds an embedded disc.

Thus there is a disc $\mathbb{D} \subset \Omega_t \subset \kappa^{-1}(t)$ about $(x^*,x^*,z^*)$, of $w(t)$-good characters, for the specified orientation. Since $\kappa^{-1}(t)$ is invariant under the reflection $(x,y,z) \mapsto (y,x,z)$, and $(x^*,x^*,z^*)$ is a fixed point, we may take $\mathbb{D}$ invariant under this reflection. Note $\mu_t(\mathbb{D})>0$: the measure $\mu_t$ is obtained by integrating the symplectic form $\omega_t$ on $X_t(S_1)$ (see \cite{Goldman03}). Hence the $\mu_t$-measure of $w(t)$-good characters in $\Omega_t$ for the specified orientation is nonzero.

The orbit of $\mathbb{D} \subset \Omega_t$ under $\Gamma$ is not null; so by ergodicity theorem \ref{thm:Goldman_ergodicity}, this orbit is conull in $\Omega_t$. Thus almost every character in $\Omega_t$ is equivalent to one in $\mathbb{D}$ by some change of basis. By post-composing with $(x,y,z) \mapsto (y,x,z)$ if necessary, which is the action of the orientation-reversing change of basis $(G,H) \mapsto (H,G)$, in fact almost every character in $\Omega_t$ is equivalent to one in $\mathbb{D}$ by an \emph{orientation-preserving} change of basis. But the definition of $w(t)$-goodness is clearly invariant under the action of any orientation-preserving change of basis. Thus almost every character in $\Omega_t$ is $w(t)$-good, and we are done.
\end{Proof}

\subsection{Piecing together character varieties}

It remains to see how the character varieties and associated
measures decompose when we cut and paste our surfaces. As our closed
surface $S$ is cut along a curve $D$ into a punctured torus $S_1$
and another surface $W$, we obtain natural maps between spaces and character varieties.
\[
    \begin{array}{ccccccccc}
        D   &   \rightarrow &   S_1         	&&&&         X(D)    &   \leftarrow  &   X(S_1) \\
        \downarrow &        &   \downarrow  	&&&&        \uparrow    &           &   \uparrow \\
        W &   \rightarrow  &   S		&&&&        X^\pm(W)  &   \leftarrow  &   X^\pm(S \; | \; D).
    \end{array}
\]
Here the pushout
\[
    X^\pm(S \; | \; D) = \left\{ \left( [\rho_1], [\rho_W] \right) \in
    X(S_1) \times X^\pm(W) \; | \; [\rho_1|_{\pi_1(D)}] =
    [\rho_W|_{\pi_1(D)}] \in X(D) \right\}
\]
is \emph{not} the same as $X^\pm(S)$; for instance, for holonomy
representations $\rho_1, \rho_2$ with the same trace along $D$,
there are many possible representations on $S$ corresponding to
twisting around the curve $D$. However there is a natural map $X^\pm(S)
\To X^\pm(S \; | \; D)$, and hence a composition $X^\pm(S) \To X^\pm(S \; | \; D) \To X^\pm(S_1)$. (Throughout we write $X^\pm$ to recall that these character varieties have coordinates given by traces of the matrix images of various curves in $SL_2\R$; the matrices corresponding to a relator multiply to $1 \in PSL_2\R$ but this may lift to $\pm 1 \in SL_2\R$, i.e. a twisted representation. Since we regard $\pi_1(S_1)$ as a free group there is no relator and need for twisted representations.)

Away from singularities, which have measure zero, the map $X^\pm(S) \To
X(S_1)$ is a submersion, since it can be taken to be a polynomial
map, indeed a coordinate map, from a $(6g-6)$-dimensional set to a
$3$-dimensional set. Recall the character variety is defined by
taking traces of a fixed set of curves on the surface $S$. As usual, we take for $S_1$ a set of standard curves $(G,H,GH)$ on
$S_1$, where $G,H$ is a basis of $\pi_1(S)$. We can take the chosen
curves on $S$ to contain the chosen curves on $S_1$ so that the map
$X^\pm(S) \To X(S_1)$ is just a coordinate projection. Let the
coordinates on $X(S_1)$ be $(x,y,z)$, and let the coordinates on
$X^\pm(S) \subset \R^{k+3}$ be $(x,y,z,w_1, \ldots, w_k)$. As $\mu_S$ is absolutely
continuous with respect to Lebesgue measure, there exists a real
function $f$ (a Radon-Nikodym derivative) such that for any Lebesgue
measurable $A \subset X^\pm(S)$, we have
\[
    \mu_S (A) = \int_{(x,y,z,w_1, \ldots, w_k) \in X^{\pm}(S)}
    \chi_A \; f \; d\lambda(x,y,z,w_1, \ldots, w_k)
\]
where $d\lambda$ denotes the $(6g-6)$-dimensional Euclidean area
form in $X^\pm(S) \subset \R^{k+3}$ and $\chi_A$ denotes the
characteristic function of the set $A$.

We claim that the symplectic 2-forms on $X^\pm(S)$ and $X(S_1)$ are
related by the natural map $X^\pm(S) \To X(S_1)$. As described in
section \ref{sec:character_variety}, the tangent space to $X^\pm(S)$ at a point $[\rho_0]$ is
$H^1(S;\B)$, where $\B$ is the bundle of coefficients over $S$
associated with the $\pi_1(S)$-module $\mathfrak{sl}_2\R_{\Ad \rho}$.
The tangent space to $X(S_1)$ is likewise $H^1(S_1; \B_1)$ where
$\B_1$ is the bundle of coefficients over $S_1$ associated with the
$\pi_1(S_1)$-module $\mathfrak{sl}_2\R_{\Ad \rho_1}$, where $\rho_1$ is
the induced homomorphism on $S_1$. Note $\B_1 = \B|_{S_1}$. So the
natural map $\iota: S_1 \hookrightarrow S$ induces $\iota^*:
H^1(S;\B) \To H^1(S_1;\B_1)$, and by naturality of cup product (see
e.g. \cite{Hodgson_thesis}) we obtain a commutative diagram
\[
    \begin{array}{ccc}
    T_{[\rho]} X(S) \times T_{[\rho]} X(S) \cong H^1(S;\B) \times
    H^1(S; \B) && \\
    \downarrow \iota^* \times \iota^* &
    \stackrel{\cup}{\searrow} & \\
    T_{[\rho_1]} X(S_1) \times T_{[\rho_1]} X(S_1) \cong H^1 (S_1;
    \B_1) \times H^1(S_1; \B_1) & \stackrel{\cup}{\rightarrow} & \R.
    \end{array}
\]

\begin{Proof}[of \ref{measure_theory_prop}]
    Recall $B = \cup_{D \in \mathcal{U}} B_D$ and $B_D = \cup_t
    B_{D,t}$. For a separating curve $D \in \mathcal{U}$, splitting $S$ into
    a punctured torus $S_1$ and a surface $W$, we defined $B_{D,t} \subset X^\pm (S)$ to be the set of all characters with euler class $-\chi(S)-1$, which take $D$ to have trace $t$, and which restrict on $S_1$ to a character in $B_t$.

    We first show $\mu_S(B_D) = 0$ for given $D \in \mathcal{U}$.

    Under the natural coordinate projection $j: X^\pm (S) \To X (S_1)$, the image $A$ of $B_{D}$ is a set of characters of representations of
    $\pi_1(S_1)$, in the various $X_t(S_1)$, which are $w(t)$-bad, for a specified orientation of $S_1$. The image of each set $B_{D,t}$ under $j$ lies in $X_t(S_1)$ and is $B_t$; thus $A = \cup_{t>2} B_t$. By by proposition \ref{ergodicity_prop} we have $\mu_t(B_t) = 0$. Note that $B_D \subseteq \left( A \times \R^k \right) \cap X^\pm(S)$. We have:
    \begin{align*}
        \mu_{S}(B_D) &= \int_{(x,y,z, w_1, \ldots, w_k) \in X^\pm(S)}
        \chi_{B_D} \; f \; d\lambda(x,y,z,w_1, \ldots, w_k) \\
        & \leq \int_{(x,y,z,w_1, \ldots, w_k) \in X^\pm(S)} \chi_{ \left( A
        \times \R^k \right) \cap X(S)} \; f \; d\lambda(x,y,z,w_1,
        \ldots, w_k) \\
        &= \int_{(w_1, \ldots, w_k)} \left( \int_{(x,y,z) \in
        X(S_1)} \chi_A \; f(x,y,z,w_1, \ldots, w_k) \; d\lambda(x,y,z)
        \right) d\lambda(w_1, \ldots, w_k).
    \end{align*}
    Thus it is sufficient to show that for any given $(w_1, \ldots,
    w_k)$, the inner integral is zero.

    Introduce the variable $t = \kappa(x,y,z) = \Tr \rho(D)$.
    The map $(x,y,z) \mapsto \kappa(x,y,z) = t$ is polynomial, hence
    measurable, so we may disintegrate the measure $d\lambda(x,y,z)$ over $t$ and
    obtain a family of measures on the level sets $X_t(S_1)$ (for details
    see e.g. \cite{Pollard}). On the level set $X_t(S_1) = \kappa^{-1}(t)$, we have the
    symplectic 2-form $\omega_t$ and measure $\mu_t$. But we have
    seen above that by naturality of the cup product, $\omega_t$ is the projection of $\omega$ under the natural map $X^\pm(S) \To
    X(S_1)$. Hence over each $X_t(S_1)$ we have an integral of $\chi_{B_t}$, times some function, with respect to $\mu_t$. Since $\mu_t(B_t) = 0$ the integral is $0$ for each $t$; integrating over $t$ we still have zero. Thus $\mu_S(B_D) = 0$.

    Now $\mathcal{U}$ certainly has cardinality no greater than the
    fundamental group of $S$, hence is countable. So the union
    $B = \cup_D B_D \subset X(S)$ is a countable union of sets of measure
    zero, hence has measure zero.
\end{Proof}

\addcontentsline{toc}{section}{References}

\small

\bibliography{danbib}
\bibliographystyle{amsplain}

\end{document}